\documentclass[12pt]{amsart}

\usepackage{enumitem}
\usepackage{psfrag}
\usepackage{graphicx}
\usepackage{pinlabel}
\usepackage{amsmath}
\usepackage{amssymb}
\usepackage{amscd}
\usepackage{picinpar}
\usepackage{times}
\usepackage{pb-diagram}
\usepackage{wrapfig}
\usepackage{xspace}
\usepackage{hyperref}
\usepackage{url}
\usepackage{caption}
\usepackage{subcaption}
\usepackage{fancyhdr}
\usepackage{dsfont}
\usepackage{mathrsfs}
\usepackage{amsfonts}
\usepackage{amsthm}
\usepackage{xcolor}
\usepackage{appendix}
\usepackage{cite}

\usepackage{tikz}

\usetikzlibrary{calc,shapes,arrows}

\usepackage{pgfplots}

\usetikzlibrary{patterns}

\usepgfplotslibrary{dateplot}

\newtheorem{theorem}{Theorem}[section]
\newtheorem{lemma}[theorem]{Lemma}

\newtheorem{corollary}[theorem]{Corollary}
\newtheorem{definition}[theorem]{Definition}
\newtheorem{conjecture}[theorem]{Conjecture}
\newtheorem{remark}[theorem]{Remark}
\newtheorem*{theorem*}{Theorem}
\usepackage{overpic}
\newcommand{\bal}{\begin{aligned}}      \newcommand{\eal}{\end{aligned}}
\newcommand{\ba}{\begin{array}}      \newcommand{\ea}{\end{array}}
\newcommand{\bc}{\begin{center}}     \newcommand{\ec}{\end{center}}
\newcommand{\be}{\begin{enumerate}}  \newcommand{\ee}{\end{enumerate}}
\newcommand{\beq}{\begin{eqnarray}}  \newcommand{\eeq}{\end{eqnarray}}
\newcommand{\beQ}{\begin{eqnarray*}} \newcommand{\eeQ}{\end{eqnarray*}}
\newcommand{\bi}{\begin{itemize}}    \newcommand{\ei}{\end{itemize}}
\newcommand{\bt}{\begin{tabular}}    \newcommand{\et}{\end{tabular}}
\newcommand{\bdm}{\begin{displaymath}} \newcommand{\edm}{\end{displaymath}}

\def\qed{\hfill{Q.E.D.}\smallskip}

\begin{document}

\title{Rigidity of discrete conformal structures on surfaces}
\author{Xu Xu}

\address{School of Mathematics and Statistics, Wuhan University, Wuhan, 430072, P.R.China}
 \email{xuxu2@whu.edu.cn}



\begin{abstract}
In \cite{G3}, Glickenstein introduced the discrete conformal structures on polyhedral surfaces in an axiomatic approach from Riemannian geometry perspective.
Glickenstein's discrete conformal structures include Thurston's circle packings, Bowers-Stephenson's inversive distance circle packings
and Luo's vertex scalings as special cases.
Glickenstein \cite{G5} further conjectured the rigidity of the discrete conformal structures on polyhedral surfaces.
Glickenstein's conjecture includes Luo's conjecture on the rigidity of vertex scalings \cite{L1} and
Bowers-Stephenson's conjecture on the rigidity of inversive distance circle packings \cite{BSt} as special cases.
In this paper, we prove Glickenstein's conjecture using variational principles.
This unifies and generalizes the well-known results of Luo \cite{L4} and Bobenko-Pinkall-Springborn \cite{BPS}.
Our method provides a unified approach to similar problems.
We further discuss the relationships of Glickenstein's discrete conformal structures on polyhedral surfaces and $3$-dimensional hyperbolic geometry.
As a result, we obtain some new results on the convexities of the
co-volume functions of some generalized $3$-dimensional hyperbolic tetrahedra.
\end{abstract}

\maketitle


\section{Introduction}
Discrete conformal structure on polyhedral manifolds is a discrete analogue of the conformal structure on Riemannian manifolds,
which assigns the discrete metrics by scalar functions defined on the vertices.
Since the work of Thurston \cite{T1}, there have been lots of researches
on different types of discrete conformal structures on polyhedral surfaces,
including the tangential circle packings, Thurston's circle packings,
Bowers-Stephenson's inversive distance circle packings, Luo's vertex scalings and others.
Most of these discrete conformal structures were invented and studied individually in the literature.
In \cite{G3}, Glickenstein developed an axiomatic approach to the Euclidean discrete conformal structures on polyhedral surfaces from Riemannian geometry perspective.
Following Glickenstein's original work \cite{G3}, Glickenstein-Thomas \cite{GT}  introduced the hyperbolic and spherical discrete conformal structures on polyhedral surfaces
in an axiomatic approach.
Glickenstein-Thomas \cite{GT} further studied the classification of Glickenstein's discrete conformal structures on polyhedral surfaces.
See also Xu-Zheng \cite{XZ4} for a complete classification of Glickenstein's discrete conformal structures.
According to the classification, Glickenstein's discrete conformal structures
include different types of circle packings and Luo's vertex scalings on polyhedral surfaces as special cases
and generalize them to a very general context.
In this paper, we study the rigidity of Glickenstein's discrete conformal structures on closed polyhedral surfaces.
In \cite{X4}, we study the deformation of Glickenstein's discrete conformal structures on surfaces.
\subsection{Polyhedral surfaces, discrete conformal structures and
the rigidity results}
Suppose $(M, \mathcal{T})$ is a connected closed triangulated surface with a triangulation $\mathcal{T}$,
which is the quotient of a finite disjoint union of triangles by identifying all the edges of triangles
in pair by homeomorphisms. We use $V, E, F$ to denote the set of vertices, unoriented edges and faces
in $\mathcal{T}$ respectively. For simplicity, we use one index to denote a vertex (such as $i\in V$),
two indices to denote an edge (such as $\{ij\}\in E$) and three indices to denote a triangle (such as $\{ijk\}\in F$).
We further use $f_i=f(i)$ for a function $f: V\rightarrow \mathbb{R}$, $f_{ij}=f(\{ij\})$ for a function $f: E\rightarrow \mathbb{R}$,
and $f_{ijk}=f(\{ijk\})$ for a function $f: F\rightarrow \mathbb{R}$ for simplicity.
Denote  the set of positive real numbers as $\mathbb{R}_{>0}$ and $|V|=N$.

\begin{definition}[\cite{L5}]\label{defn of polyhedral surface}
A polyhedral surface $(M, \mathcal{T}, l)$ with background geometry $\mathbb{G}$ ($\mathbb{G}=\mathbb{E}^2, \mathbb{H}^2$ or $\mathbb{S}^2$)
is a triangulated surface $(M, \mathcal{T})$ with a map $l: E\rightarrow \mathbb{R}_{>0}$ such that
any face $\{ijk\}\in F$ can be embedded in $\mathbb{G}$ as a nondegenerate triangle
with edge lengths $l_{ij}, l_{ik}, l_{jk}$ given by $l$.
We call $l: E\rightarrow \mathbb{R}_{>0}$ as a Euclidean (hyperbolic or spherical respectively) polyhedral metric
if $\mathbb{G}=\mathbb{E}^2$ ($\mathbb{G}=\mathbb{H}^2$ or $\mathbb{G}=\mathbb{S}^2$ respectively).
\end{definition}

The nondegenerate condition for the face $\{ijk\}\in F$ in Definition \ref{defn of polyhedral surface}
is equivalent to the edge lengths $l_{ij}, l_{ik}, l_{jk}$ satisfy the triangle inequalities
($l_{ij}+l_{ik}+l_{jk}<2\pi$ additionally if $\mathbb{G}=\mathbb{S}^2$).
Intuitively, a polyhedral surface with background geometry $\mathbb{G}$ ($\mathbb{G}=\mathbb{E}^2, \mathbb{H}^2$ or $\mathbb{S}^2$) can be obtained
by gluing triangles in $\mathbb{G}$ isometrically along the edges in pair.
For polyhedral surfaces, there may exist conic singularities at the vertices,
which can be described by combinatorial curvatures.
The combinatorial curvature is a map $K: V\rightarrow (-\infty, 2\pi)$ that assigns the vertex $i\in V$
$2\pi$ less the sum of inner angles at $i$, i.e.
\begin{equation}
K_i=2\pi-\sum_{\{ijk\}\in F}\theta_{i}^{jk},
\end{equation}
where $\theta_{i}^{jk}$ is the inner angle at $i$ in the triangle $\{ijk\}$.

\begin{definition}[\cite{G3, GT}]\label{defn of discrete conformal structure}
Suppose $(M, \mathcal{T})$ is a triangulated connected closed surface and $\varepsilon: V\rightarrow \{-1, 0, 1\}$,
$\eta: E\rightarrow \mathbb{R}$ are two weights defined on the vertices and edges respectively.
A discrete conformal structure on the weighted triangulated surface $(M, \mathcal{T}, \varepsilon, \eta)$ with background geometry $\mathbb{G}$
is composed of  the maps $f: V\rightarrow \mathbb{R}$ such that
\begin{description}
  \item[(1)] the edge length $l_{ij}$ for the edge $\{ij\}\in E$ is given by
  \begin{equation}\label{defn of Euclidean length}
  \begin{aligned}
  l_{ij}=\sqrt{\varepsilon_i e^{2f_i}+\varepsilon_j e^{2f_j}+2\eta_{ij}e^{f_i+f_j}}
  \end{aligned}
  \end{equation}
  for $\mathbb{G}=\mathbb{E}^2$,
  \begin{equation}\label{defn of hyperbolic length}
  \begin{aligned}
  l_{ij}=\cosh^{-1}\left(\sqrt{(1+\varepsilon_ie^{2f_i})(1+\varepsilon_je^{2f_j})}+\eta_{ij}e^{f_i+f_j}\right)
  \end{aligned}
  \end{equation}
  for $\mathbb{G}=\mathbb{H}^2$ and
  \begin{equation}\label{defn of spherical length}
  \begin{aligned}
  l_{ij}=\cos^{-1}\left(\sqrt{(1-\varepsilon_ie^{2f_i})(1-\varepsilon_je^{2f_j})}-\eta_{ij}e^{f_i+f_j}\right)
  \end{aligned}
  \end{equation}
  for $\mathbb{G}=\mathbb{S}^2$;
  \item[(2)]  the edge length function  $l: E\rightarrow \mathbb{R}_{>0}$ defined by
  (\ref{defn of Euclidean length}), (\ref{defn of hyperbolic length}), (\ref{defn of spherical length})
  is a Euclidean, hyperbolic and spherical polyhedral metric on $(M, \mathcal{T})$ respectively.
\end{description}
The weights $\varepsilon: V\rightarrow \{-1, 0, 1\}$ and $\eta: E\rightarrow \mathbb{R}$ are called the scheme coefficient and discrete conformal coefficient respectively.
A function $f: V\rightarrow \mathbb{R}$ is called a discrete conformal factor
and a function $f: V\rightarrow \mathbb{R}$ with the induced edge length function $l: E\rightarrow \mathbb{R}_{>0}$ being
a polyhedral metric is called a nondegenerate discrete conformal factor.
\end{definition}


\begin{remark}
It is a remarkable result of Glickenstein-Thomas \cite{GT} that Glickenstein's discrete conformal structure  can be classified,
which has the form
given in Definition \ref{defn of discrete conformal structure} with $\varepsilon_i\in \{-1,0,1\}$ replaced by a constant $\alpha_i\in \mathbb{R}$.
See also Xu-Zheng \cite{XZ4}.
As pointed out by Thomas (\cite{T} page 53), one can reparameterize Glickenstein's discrete conformal structures so that $\alpha_i\in \{-1, 0, 1\}$
while keeping the induced polyhedral metrics invariant.
This is the motivation of Definition \ref{defn of discrete conformal structure}.
\end{remark}
\begin{table}[ht]
\centering
\begin{tabular}{|l l l l|}
\hline
  Scheme & $\varepsilon_i$ & $\varepsilon_j$ & $\eta_{ij}$ \\ \hline
  Tangential circle packings & $+1$ & $+1$ & $+1$ \\
  Thurston's circle packings & $+1$ & $+1$ & $(-1,1]$ \\
  Bowers-Stephenson's inversive distance circle packings\ \  & $+1$ & $+1$ & $(-1, +\infty)$ \\
  Luo's vertex scalings & $0$ & $0$ & $(0, +\infty)$ \\
  Glickenstein's discrete conformal structures & $\{+1, 0, -1\}$ & $\{+1, 0, -1\}$ & $(-1, +\infty)$ \\
  \hline
\end{tabular}
\caption{Different types of discrete conformal structures} 
\label{Relationships of DCS} 
\end{table}
\begin{remark}

The relationships of Glickenstein's discrete conformal structures in Definition \ref{defn of discrete conformal structure}
and the existing special types of discrete conformal structures are contained in Table \ref{Relationships of DCS}.
By Table \ref{Relationships of DCS}, the tangential circle packing is a special case of Thurston's circle packing and
Thurston's circle packing is a special case of Bowers-Stephenson's inversive distance circle packing.
For simplicity, we unify all these three types of circle packings as inversive distance circle packings in the following.
By Table \ref{Relationships of DCS} again, Glickenstein's discrete conformal structures in Definition \ref{defn of discrete conformal structure} include
Bowers-Stephenson's inversive distance circle packings and Luo's vertex scalings as special cases.
Furthermore, Glickenstein's discrete conformal structures in Definition \ref{defn of discrete conformal structure}
include the mixed type discrete conformal structures. Specially, it contains the type with
$\varepsilon_i=0$ for some vertices $i\in V$ and $\varepsilon_j=1$ for the other vertices $j\in V$.
The geometry of such mixed type discrete conformal structures is seldom studied in the literature.

\end{remark}

A basic problem in discrete conformal geometry is to understand the relationships between the discrete conformal factors
and their combinatorial curvatures.
We prove the following result on the rigidity of Glickenstein's discrete conformal structures on polyhedral surfaces.
\begin{theorem}\label{main rigidity introduction}
Suppose $(M, \mathcal{T}, \varepsilon, \eta)$ is a weighted triangulated connected closed surface with
the weights $\varepsilon: V\rightarrow \{0, 1\}$ and $\eta: E\rightarrow \mathbb{R}$ satisfying the structure conditions
\begin{equation}\label{structure condition 1}
\begin{aligned}
\varepsilon_s \varepsilon_t +\eta_{st}>0, \ \ \forall \{st\}\in E
\end{aligned}
\end{equation}
and
\begin{equation}\label{structure condition 2}
\begin{aligned}
\varepsilon_q\eta_{st}+\eta_{qs}\eta_{qt}\geq 0, \ \ \{q, s, t\}= \{i, j, k\}
\end{aligned}
\end{equation}
for any triangle $\{ijk\}\in F$.
\begin{description}
  \item[(a)] A  nondegenerate Euclidean discrete conformal factor $f: V\rightarrow \mathbb{R}$ on $(M, \mathcal{T}, \varepsilon, \eta)$ is determined by its combinatorial
  curvature $K: V\rightarrow \mathbb{R}$ up to a vector $c(1,1, \cdots, 1), c\in \mathbb{R}$.
  \item[(b)] A nondegenerate hyperbolic discrete conformal factor $f: V\rightarrow \mathbb{R}$ on $(M, \mathcal{T}, \varepsilon, \eta)$ is determined by its combinatorial
  curvature $K: V\rightarrow \mathbb{R}$.
\end{description}
\end{theorem}

\begin{remark}
Theorem \ref{main rigidity introduction} confirms a conjecture of Glickenstein in \cite{G5}.
If $\varepsilon\equiv1$, Theorem \ref{main rigidity introduction} is reduced to the rigidity of
Bowers-Stephenson's inversive distance circle packings on surfaces
obtained by Guo \cite{Guo}, Luo \cite{L4} and the author \cite{X1,X3},
which was conjectured by Bowers-Stephenson \cite{BSt}.
If $\varepsilon\equiv0$, Theorem \ref{main rigidity introduction} is reduced to the rigidity of
Luo's vertex scalings on surfaces obtained by Luo \cite{L1} and Bobenko-Pinkall-Springborn \cite{BPS},
the global rigidity of which was conjectured by Luo \cite{L1}.
Theorem \ref{main rigidity introduction} unifies these rigidity results.
Furthermore, Theorem \ref{main rigidity introduction} includes the rigidity of the mixed type
discrete conformal structures, for which $\varepsilon_i=1$ for some vertices $i\in V_1\neq \emptyset$
and $\varepsilon_j=0$ for the other vertices $j\in V\setminus V_1\neq \emptyset$.
The local rigidity of Glickenstein's discrete conformal structures on polyhedral surfaces was first proved by
Glickenstein \cite{G3} and Glickenstein-Thomas \cite{GT}
under a  condition that the discrete conformal structures induce a well-centered geometric center for each triangle
in the triangulation.
The local rigidity for some subcases of Glickenstein's discrete conformal structures was also proved by
Guo-Luo \cite{GL}. Theorem \ref{main rigidity introduction} includes these results on local rigidity as special cases.
\end{remark}

\subsection{Relationships with $3$-dimensional hyperbolic geometry}

Motivated by Bobenko-Pinkall-Springborn's observations \cite{BPS} on the deep relationships of Luo's vertex scalings
on polyhedral surfaces and $3$-dimensional hyperbolic geometry,
Zhang-Guo-Zeng-Luo-Yau-Gu \cite{ZGZLYG} constructed
Glickenstein's discrete conformal structures via generalized $3$-dimensional hyperbolic tetrahedra.
The basic strategy is to construct a generalized hyperbolic tetrahedron $T_{Oijk}$ with the vertices $O, v_i, v_j, v_k$ in $\mathbb{H}^3$, ideal or hyper-ideal. And then the discrete conformality naturally appears at some vertex triangle.
In this paper, we focus on the case that the vertices are ideal or hyper-ideal.
The vertex $O$ is ideal when we study Glickenstein's Euclidean discrete conformal structures, and hyper-ideal when we study Glickenstein's
hyperbolic discrete conformal structures.
The vertex $v_s\in \{v_i, v_j, v_k\}$ is hyper-ideal if $\varepsilon_s=1$, and ideal if $\varepsilon_s=0$.
In the case that $v_s\in \{v_i, v_j, v_k\}$ is hyper-ideal, the line segment $Ov_s$ is required to have nonempty intersection with $\mathbb{H}^3$ in the Klein model.
For each pair $\{v_s, v_t\}\subseteq\{v_i, v_j, v_k\}$, a weight $\eta_{st}$ can be naturally assigned via the signed edge length of $\{v_sv_t\}$.
In the Euclidean background geometry, the edge lengths of the vertex triangle $T_{Oijk}\cap H_O$ is given by  (\ref{defn of Euclidean length}),
where $H_O$ is the horosphere attached to the ideal vertex $O$ and $f_s$ is minus of the signed decorated edge length $l_{Ov_s}$ with $s\in \{i,j,k\}$.
The case for hyperbolic background geometry is similar.
By truncating the generalized hyperbolic tetrahedron $T_{Oijk}$ with hyperbolic planes dual to the hyper-ideal vertices,
we can attach it with a generalized hyperbolic polyhedron $P$ with finite volume.
For the details on the construction of $T_{Oijk}$ and assignments of $\eta_{ij}, \eta_{ik}, \eta_{jk}$,
please refer to Section \ref{section 4}.


\begin{theorem}\label{theorem convexity of co-volume}
Suppose $T=\{Oijk\}$ is a generalized tetrahedron constructed above.
\begin{description}
  \item[(a)] The weights $\eta_{ij}, \eta_{ik}, \eta_{jk}$ on the edges $\{ij\}, \{ik\}, \{jk\}$
satisfy the structure conditions (\ref{structure condition 1}) and  (\ref{structure condition 2}).
  \item[(b)] The co-volume of the generalized tetrahedron $T=\{Oijk\}$
  defined by (\ref{definition of covolume})
  with fixed weights $\eta_{ij}, \eta_{ik}, \eta_{jk}$
  is a convex function of the signed edge lengths $l_{Ov_i}, l_{Ov_j}, l_{Ov_k}$.
\end{description}
\end{theorem}

\subsection{Basic ideas of the proof of Theorem \ref{main rigidity introduction}}

The proof for the rigidity of Glickenstein's discrete conformal structures on triangulated surfaces, i.e. Theorem \ref{main rigidity introduction},
involves a variational principle
introduced by Colin de Verdi\`{e}re \cite{DV}.
The variational principle  has been extensively studied in
\cite{BPS, BSp, CL,DGL,Guo, GL, Lei, L1,L4,L5, R1,S} and others.
Glickenstein \cite{G3} and Glickenstein-Thomas \cite{GT} generalized the variational principle to
Glickenstein's discrete conformal structures and proved
some results on the local rigidity of Glickenstein's discrete conformal structures.
In this paper, we use Glickenstein's variational principle to prove the local and global rigidity of Glickenstein's discrete conformal structures.
The key ingredient using Glickenstein's variational principle to prove the rigidity
is constructing a globally defined convex function with the combinatorial curvature as its gradient.
This can be reduced to constructing a globally defined concave function of the discrete conformal factors
on a triangle with the inner angles as its gradient.
The main difficulties come from the characterization of the admissible space of nondegenerate discrete conformal factors on a triangle and
the local concavity of the function with the inner angles as its gradient.
We construct such a function in three steps.
In the first step, we give an analytical characterization of the admissible space of nondegenerate discrete conformal factors on a triangle.
This is accomplished by solving the global triangle inequalities
with the help of the geometric center introduced by Glickenstein \cite{G3,G4}.
As a result, the admissible space of nondegenerate discrete conformal factors on a triangle is proved to be homotopy equivalent to $\mathbb{R}^3$ and hence simply connected.
This implies that the Ricci energy function, defined as the integral of the inner angles on the admissible space of nondegenerate discrete conformal factors for a triangle, is well-defined.
In the second step, we show that the Ricci energy function for a triangle is locally concave.
To achieve this, we introduce the parameterized admissible space of the nondegenerate discrete conformal factors,
and choose some good point in the space such that the hession matrix of the Ricci energy function is negative definite at this point.
By the continuity of the eigenvalues of the hession matrix, we prove the local concavity of the Ricci energy function.
In the final step, we extend the locally concave Ricci energy function defined on the admissible space of nondegenerate discrete conformal factors
 for a triangle to be a globally defined concave function.
The extension is now standard since Bobenko-Pinkall-Springborn's important work \cite{BPS}.
In this paper, we use Luo's generalization \cite{L4} of Bobenko-Pinkall-Springborn's extension in \cite{BPS} to
extend the locally concave Ricci energy function to be a globally defined concave function.


\subsection{The organization of the paper}
In Section \ref{section 2}, we study the rigidity of Glickenstein's Euclidean discrete conformal structures on polyhedral surfaces and
prove a generalization of Theorem \ref{main rigidity introduction} (a).
In Section \ref{section 3}, we study the rigidity of Glickenstein's  hyperbolic discrete conformal structures on polyhedral surfaces and
prove a generalization of Theorem \ref{main rigidity introduction} (b).
In Section \ref{section 4}, we discuss the relationships of Glickenstein's discrete conformal structures on polyhedral surfaces
and $3$-dimensional hyperbolic geometry and prove Theorem \ref{theorem convexity of co-volume}.
In Section \ref{section 5}, we discuss some open problems.
\\
\\
\textbf{Acknowledgements}\\[8pt]
The research of the author is supported by Fundamental Research Funds for the Central Universities under Grant no. 2042020kf0199.

\section{Euclidean discrete conformal structures}\label{section 2}
\subsection{Admissible space of nondegenerate Euclidean discrete conformal factors on a triangle}\label{subsection Euclidean admissible space}
Let $\sigma=\{ijk\}$ be a triangle with the vertex set $V_\sigma=\{i,j,k\}$ and the edge set $E_\sigma=\{\{ij\}, \{ik\},\{jk\}\}$.
Unless otherwise declared,
we use $(\{ijk\}, \varepsilon,\eta)$ to denote a weighted triangle with two weights $\varepsilon: V_\sigma\rightarrow \{0, 1\}$ and $E_\sigma\rightarrow \mathbb{R}$
satisfying the structure conditions (\ref{structure condition 1}) and (\ref{structure condition 2}) in the following.
In the Euclidean background geometry,
the lengths $l_{ij}, l_{ik}, l_{jk}$ of the edges in $E_\sigma$ are defined by the discrete conformal factor $f: V_\sigma \rightarrow \mathbb{R}$ via the formula (\ref{defn of Euclidean length}).
The discrete conformal factor $f: V_\sigma \rightarrow \mathbb{R}$  is \textit{nondegenerate} if $l_{ij}, l_{ik}, l_{jk}$ satisfy the triangle inequalities, otherwise it is \textit{degenerate}.
We use $\Omega_{ijk}^E(\eta)$ to denote the space of nondegenerate Euclidean discrete conformal factors for the weighted triangle $(\{ijk\}, \varepsilon,\eta)$.

Due to $\varepsilon\in \{0,1\}$ and the structure condition (\ref{structure condition 1}),
the Cauchy inequality implies
$\varepsilon_i e^{2f_i}+\varepsilon_j e^{2f_j}+2\eta_{ij}e^{f_i+f_j}\geq 2(\varepsilon_i\varepsilon_j+\eta_{ij})e^{f_i+f_j}>0$.
Therefore, the Euclidean edge length $l_{ij}$ in (\ref{defn of Euclidean length}) is well-defined.
Note that the edge lengths $l_{ij}, l_{ik}, l_{jk}$ satisfy
the triangle inequalities
\begin{equation}\label{Euclidean triangle inequality}
\begin{aligned}
l_{ij}<l_{ik}+l_{jk},\   l_{ik}<l_{ij}+l_{jk},\  l_{jk}<l_{ij}+l_{ik}
\end{aligned}
\end{equation}
if and only if
\begin{equation}\label{equivalent Euclidean triangle inequality}
\begin{aligned}
0<&(l_{ij}+l_{ik}+l_{jk})(l_{ij}+l_{ik}-l_{jk})(l_{ij}-l_{ik}+l_{jk})(-l_{ij}+l_{ik}+l_{jk})\\
=&2l_{ij}^2l_{ik}^2+2l_{ij}^2l_{jk}^2+2l_{ik}^2l_{jk}^2-l_{ij}^4-l_{ik}^4-l_{jk}^4.
\end{aligned}
\end{equation}
For simplicity, set
\begin{equation}\label{notation ri and fi in the Euclidean case}
\begin{aligned}
r_i=e^{f_i}, \ \forall i \in V_\sigma.
\end{aligned}
\end{equation}
Then the edge length $l_{ij}$ in the Euclidean background geometry is given by
\begin{equation}\label{definition of Euclidean length}
\begin{aligned}
l_{ij}=\sqrt{\varepsilon_ir_i^2+\varepsilon_jr_j^2+2\eta_{ij}r_ir_j}.
\end{aligned}
\end{equation}
The vector $r=(r_i,r_j,r_k)\in \mathbb{R}^{3}_{>0}$ is called as a \textit{radius vector}.
Paralleling to the discrete conformal factors,
a radius vector $r: V_\sigma \rightarrow \mathbb{R}_{>0}$ is \textit{nondegenerate} if $l_{ij}, l_{ik}, l_{jk}$ satisfy the triangle inequalities, otherwise it is \textit{degenerate}.
Submitting (\ref{definition of Euclidean length}) into (\ref{equivalent Euclidean triangle inequality}) and
by direct calculations, we have
\begin{equation}\label{Euclidean triangle inequa comp}
\begin{aligned}
&(l_{ij}+l_{ik}+l_{jk})(l_{ij}+l_{ik}-l_{jk})(l_{ij}-l_{ik}+l_{jk})(-l_{ij}+l_{ik}+l_{jk})\\
=&4r_i^2r_j^2r_k^2[(\varepsilon_i\varepsilon_j-\eta_{ij}^2)r_k^{-2}+(\varepsilon_i\varepsilon_k-\eta_{ik}^2)r_j^{-2}+(\varepsilon_j\varepsilon_k-\eta_{jk}^2)r_i^{-2}\\
  &+2(\varepsilon_k\eta_{ij}+\eta_{ik}\eta_{jk})r_i^{-1}r_j^{-1}+2(\varepsilon_j\eta_{ik}+\eta_{ij}\eta_{jk})r_i^{-1}r_k^{-1}
  +2(\varepsilon_i\eta_{jk}+\eta_{ij}\eta_{ik})r_j^{-1}r_k^{-1}].
\end{aligned}
\end{equation}
Set
\begin{equation}\label{kappa}
\begin{aligned}
&\kappa_i=r_i^{-1}, \kappa_j=r_j^{-1}, \kappa_k=r_k^{-1},
\end{aligned}
\end{equation}
\begin{equation}\label{gamma}
\begin{aligned}
\gamma_i=\varepsilon_{i}\eta_{jk}+\eta_{ij}\eta_{ik}, \gamma_j=\varepsilon_{j}\eta_{ik}+\eta_{ij}\eta_{jk},
\gamma_k=\varepsilon_{k}\eta_{ij}+\eta_{ik}\eta_{jk},
\end{aligned}
\end{equation}
\begin{equation}\label{Q}
\begin{aligned}
Q^E=&(\varepsilon_j\varepsilon_k-\eta_{jk}^2)\kappa_i^2+(\varepsilon_i\varepsilon_k-\eta_{ik}^2)\kappa_j^2+(\varepsilon_i\varepsilon_j-\eta_{ij}^2)\kappa_k^{2}
   +2\kappa_i\kappa_j\gamma_k+2\kappa_i\kappa_k\gamma_j+2\kappa_j\kappa_k\gamma_i.
\end{aligned}
\end{equation}
Then the structure condition (\ref{structure condition 2}) for the triangle $\sigma=\{ijk\}$ is equivalent to
\begin{equation}\label{gamma sign}
\begin{aligned}
\gamma_i\geq 0, \gamma_j\geq 0, \gamma_k\geq 0
\end{aligned}
\end{equation}
and the formula (\ref{Euclidean triangle inequa comp}) can be written as
\begin{equation*}
\begin{aligned}
(l_{ij}+l_{ik}+l_{jk})(l_{ij}+l_{ik}-l_{jk})(l_{ij}-l_{ik}+l_{jk})(-l_{ij}+l_{ik}+l_{jk})=4r_i^2r_j^2r_k^2Q^E.
\end{aligned}
\end{equation*}
As a consequence of the arguments above, we have the following result.
\begin{lemma}\label{lemma triangle equality euqi to Q>0}
For the weighted triangle $(\{ijk\}, \varepsilon,\eta)$, the edge lengths
$l_{ij}, l_{ik}, l_{jk}$ defined by (\ref{defn of Euclidean length})
satisfy the triangle inequalities if and only if $Q^E>0$.
\end{lemma}


Set
\begin{equation}\label{hi hj hk}
\begin{aligned}
h_i=&(\varepsilon_j\varepsilon_k-\eta_{jk}^2)\kappa_i+\kappa_j\gamma_k+\kappa_k\gamma_j,\\
h_j=&(\varepsilon_i\varepsilon_k-\eta_{ik}^2)\kappa_j+\kappa_i\gamma_k+\kappa_k\gamma_i,\\
h_k=&(\varepsilon_i\varepsilon_j-\eta_{ij}^2)\kappa_k+\kappa_i\gamma_j+\kappa_j\gamma_i.
\end{aligned}
\end{equation}
Then we have
$$
Q^E=\kappa_ih_i+\kappa_jh_j+\kappa_kh_k.
$$
By Lemma \ref{lemma triangle equality euqi to Q>0},
$r=(r_i,r_j,r_k)\in \mathbb{R}^3_{>0}$ is a degenerate radius vector for the triangle $\{ijk\}$
if and only if $Q^E\leq 0$. This implies that if $r=(r_i,r_j,r_k)\in \mathbb{R}^3_{>0}$ is a degenerate radius vector, then
at least one of $h_i, h_j, h_k$ is nonpositive.
Furthermore, we have the following result on the signs of $h_i, h_j, h_k$.

\begin{lemma}\label{sign of hi hj hk}
If $r=(r_i,r_j,r_k)\in \mathbb{R}^3_{>0}$ is a degenerate radius vector on the weighted triangle $(\{ijk\}, \varepsilon,\eta)$,
then one of $h_i, h_j, h_k$ is negative and the other two are positive.
\end{lemma}
\proof
We separate the proof into two steps.

\textbf{Step 1:}
For any $r=(r_i,r_j,r_k)\in \mathbb{R}^3_{>0}$,
there is no subset $\{s,t\}\subset \{i,j,k\}$ such that $h_s\leq 0$ and $h_t\leq 0$.

Suppose otherwise $h_i\leq 0$, $h_j\leq 0$. Then by the definition of $h_i, h_j$ in (\ref{hi hj hk}), we have
\begin{equation}\label{proof lemma eq1}
\begin{aligned}
\kappa_j\gamma_k+\kappa_k\gamma_j\leq (\eta_{jk}^2-\varepsilon_j\varepsilon_k)\kappa_i,
\end{aligned}
\end{equation}
\begin{equation}\label{proof lemma eq2}
\begin{aligned}
\kappa_i\gamma_k+\kappa_k\gamma_i\leq (\eta_{ik}^2-\varepsilon_i\varepsilon_k)\kappa_j.
\end{aligned}
\end{equation}
By the structure conditions (\ref{structure condition 2}), the inequalities (\ref{proof lemma eq1}) and (\ref{proof lemma eq2}) imply
$\eta_{jk}^2-\varepsilon_j\varepsilon_k\geq 0$ and $\eta_{ik}^2-\varepsilon_i\varepsilon_k\geq 0$.
By the structure conditions (\ref{structure condition 1}) and $\varepsilon_i, \varepsilon_j, \varepsilon_k\in \{0,1\}$, this implies
\begin{equation}\label{equation 1 proof hi hj nonpositive}
\begin{aligned}
\eta_{jk}-\varepsilon_j\varepsilon_k\geq 0,\ \eta_{ik}-\varepsilon_i\varepsilon_k\geq 0.
\end{aligned}
\end{equation}
Multiplying (\ref{proof lemma eq1}) and (\ref{proof lemma eq2}) gives
$(\kappa_j\gamma_k+\kappa_k\gamma_j)(\kappa_i\gamma_k+\kappa_k\gamma_i)
\leq \kappa_i\kappa_j(\eta_{jk}^2-\varepsilon_j\varepsilon_k)(\eta_{ik}^2-\varepsilon_i\varepsilon_k).$
By the structure condition (\ref{structure condition 2}), this implies
$\kappa_i\kappa_j\gamma_k^2
\leq \kappa_i\kappa_j(\eta_{jk}^2-\varepsilon_j\varepsilon_k)(\eta_{ik}^2-\varepsilon_i\varepsilon_k),$
which is equivalent to
\begin{equation}\label{F leq 0 a}
\begin{aligned}
\gamma_k^2-(\eta_{jk}^2-\varepsilon_j\varepsilon_k)(\eta_{ik}^2-\varepsilon_i\varepsilon_k)
=
-\varepsilon_i\varepsilon_j\varepsilon_k+\varepsilon_k\eta_{ij}^2+\varepsilon_j\varepsilon_k\eta_{ik}^2
+\varepsilon_i\varepsilon_k\eta_{jk}^2+2\varepsilon_k\eta_{ij}\eta_{ik}\eta_{jk}\leq 0.
\end{aligned}
\end{equation}
Set
$$F=-\varepsilon_i\varepsilon_j\varepsilon_k+\varepsilon_k\eta_{ij}^2+\varepsilon_j\varepsilon_k\eta_{ik}^2
+\varepsilon_i\varepsilon_k\eta_{jk}^2+2\varepsilon_k\eta_{ij}\eta_{ik}\eta_{jk}.$$
Then $F\leq 0$ by (\ref{F leq 0 a}).
On the other hand, by (\ref{equation 1 proof hi hj nonpositive}), $\varepsilon_i, \varepsilon_j, \varepsilon_k\in \{0,1\}$,
 (\ref{structure condition 1}) and (\ref{structure condition 2}), we have
\begin{equation}\label{F geq 0}
\begin{aligned}
F
=&\varepsilon_k(\varepsilon_i\eta_{jk}-\varepsilon_j\eta_{ik})^2
  +(\varepsilon_i\varepsilon_j+\eta_{ij})[2\varepsilon_k\eta_{ik}\eta_{jk}+\varepsilon_k(\eta_{ij}-\varepsilon_i\varepsilon_j)]\\
=&\varepsilon_k(\varepsilon_i\eta_{jk}-\varepsilon_j\eta_{ik})^2
  +(\varepsilon_i\varepsilon_j+\eta_{ij})[\varepsilon_k(\varepsilon_k\eta_{ij}+\eta_{ik}\eta_{jk})+\varepsilon_k(\eta_{ik}\eta_{jk}
  -\varepsilon_i\varepsilon_j)]\\
=&\varepsilon_k(\varepsilon_i\eta_{jk}-\varepsilon_j\eta_{ik})^2
  +(\varepsilon_i\varepsilon_j+\eta_{ij})[\varepsilon_k\gamma_k
   +\varepsilon_k(\eta_{ik}-\varepsilon_i\varepsilon_k)(\eta_{jk}-\varepsilon_j\varepsilon_k)\\
   &\ \ \ \ \ \ \ \ \ \ \ \ \ \ \ \ \ \ \ \ \ \ \ \ \ \ \ \ \ \ \ \ \ \ \ \ \ \ \ \ \ \ \ \ \ \ \ \ \ \ \ \ \ \ +\varepsilon_i\varepsilon_k(\eta_{jk}-\varepsilon_j\varepsilon_k)
   +\varepsilon_j\varepsilon_k(\eta_{ik}-\varepsilon_i\varepsilon_k)]\\
\geq& 0.
\end{aligned}
\end{equation}
Therefore, $F=0$.

In the case of $\varepsilon_k=0$, by the structure condition (\ref{structure condition 1}), we have
\begin{equation}\label{proof equ 1}
\begin{aligned}
\varepsilon_i\varepsilon_j+\eta_{ij}>0, \eta_{ik}>0, \eta_{jk}>0.
\end{aligned}
\end{equation}
By $h_i\leq 0$, $h_j\leq 0$ and $\varepsilon_k=0$, we have
$\eta_{ik}\eta_{jk}\kappa_j+\gamma_j\kappa_k\leq \eta_{jk}^2\kappa_i$ and
$\eta_{ik}\eta_{jk}\kappa_i+\gamma_i\kappa_k\leq \eta_{ik}^2\kappa_j$.
Multiplying both sides of these two inequalities gives
$$\eta_{ik}^2\eta_{jk}^2\kappa_i\kappa_j\leq (\eta_{ik}\eta_{jk}\kappa_j+\gamma_j\kappa_k)(\eta_{ik}\eta_{jk}\kappa_i+\gamma_i\kappa_k)\leq \eta_{ik}^2\eta_{jk}^2\kappa_i\kappa_j.$$
This implies
\begin{equation}\label{proof equ 2}
\begin{aligned}
\gamma_i=\varepsilon_i\eta_{jk}+\eta_{ij}\eta_{ik}=0,\ \gamma_j=\varepsilon_j\eta_{ik}+\eta_{ij}\eta_{jk}=0,
\end{aligned}
\end{equation}
which implies $\varepsilon_i\varepsilon_j-\eta_{ij}^2=0$ by (\ref{proof equ 1}).
Note that
$\varepsilon_i\varepsilon_j-\eta_{ij}^2=(\varepsilon_i\varepsilon_j-\eta_{ij})(\varepsilon_i\varepsilon_j+\eta_{ij})$.
By $\varepsilon_i\varepsilon_j+\eta_{ij}>0$ in (\ref{proof equ 1}),
we have $\varepsilon_i\varepsilon_j-\eta_{ij}=0$. By $\varepsilon_i\varepsilon_j+\eta_{ij}>0$ in (\ref{proof equ 1}) again,
we have $\varepsilon_i\varepsilon_j>0$. Therefore, $\varepsilon_i=\varepsilon_j=\eta_{ij}=1$.
Combining this with (\ref{proof equ 2}) gives $\eta_{ik}+\eta_{jk}=0$. This contradicts
$\eta_{ik}>0, \eta_{jk}>0$ in (\ref{proof equ 1}).

In the case of $\varepsilon_k=1$, by $F=0$ and (\ref{F geq 0}), we have
  \begin{equation}\label{equation 2 proof hi hj nonpositive}
\begin{aligned}
\varepsilon_i\eta_{jk}-\varepsilon_j\eta_{ik}=\eta_{ij}+\eta_{ik}\eta_{jk}=(\eta_{ik}-\varepsilon_i)(\eta_{jk}-\varepsilon_j)=
\varepsilon_i(\eta_{jk}-\varepsilon_j)=
   \varepsilon_j(\eta_{ik}-\varepsilon_i)=0.
\end{aligned}
\end{equation}
This implies $\eta_{ik}=\varepsilon_i$ or $\eta_{jk}=\varepsilon_j$.
By $\varepsilon_k=1$ and the structure condition (\ref{structure condition 1}), we have
\begin{equation}\label{proof equ 3}
\begin{aligned}
\varepsilon_i\varepsilon_j+\eta_{ij}>0, \varepsilon_i+\eta_{ik}>0, \varepsilon_j+\eta_{jk}>0.
\end{aligned}
\end{equation}
If $\eta_{ik}=\varepsilon_i$, (\ref{proof equ 3}) implies $\varepsilon_i=\eta_{ik}=1$.
Submitting this into (\ref{equation 2 proof hi hj nonpositive}) gives $\eta_{jk}=\varepsilon_j$.
By (\ref{proof equ 3}) again, we have $\varepsilon_j=\eta_{jk}=1$.
Combining $\varepsilon_i=\eta_{ik}=1$, $\varepsilon_j=\eta_{jk}=1$ and  (\ref{equation 2 proof hi hj nonpositive}), we have
$\eta_{ij}+\varepsilon_i\varepsilon_j=\eta_{ij}+\eta_{ik}\eta_{jk}=0$.
This contradicts (\ref{proof equ 3}). The same arguments also apply to the case $\eta_{jk}=\varepsilon_j$.

Therefore, there exists no subset $\{s, t\}\subset \{i,j,k\}$ such that $h_s\leq 0$ and $h_t\leq 0$.

\textbf{Step 2:} If $r=(r_i,r_j,r_k)\in \mathbb{R}^3_{>0}$ is a degenerate radius vector for the weighted triangle $(\{ijk\}, \varepsilon,\eta)$,
then one of $h_i, h_j, h_k$ is negative and the other two are positive.

By Lemma \ref{lemma triangle equality euqi to Q>0}, if $r=(r_i,r_j,r_k)\in \mathbb{R}^3_{>0}$ is a degenerate radius vector,
we have $Q^E=\kappa_ih_i+\kappa_jh_j+\kappa_kh_k\leq 0$.
This implies at least one of $h_i, h_j, h_k$ is nonpositive.
Without loss of generality, assume $h_i\leq 0$. Then the result in step 1 implies that $h_j>0, h_k>0$.
If $h_i=0$,  $h_j>0, h_k>0$, we have $Q^E=\kappa_ih_i+\kappa_jh_j+\kappa_kh_k> 0$.
This contradicts $Q^E\leq 0$.
Therefore, $h_i<0$, $h_j>0$, $h_k>0$.
\qed

%

\begin{remark}\label{power distance}
Lemma \ref{sign of hi hj hk} has an interesting geometric explanation as follows.
For a nondegenerate radius vector $r\in \mathbb{R}^3_{>0}$ on the weighted triangle $(\{ijk\}, \varepsilon,\eta)$,
there exists a geometric center $C_{ijk}$ for the triangle $\{ijk\}$ (\cite{G4}, Proposition 4),
which has the same power distance to the vertices $i, j, k$.
  \begin{figure}[!htb]
\centering
  \includegraphics[height=.4\textwidth,width=0.56\textwidth]{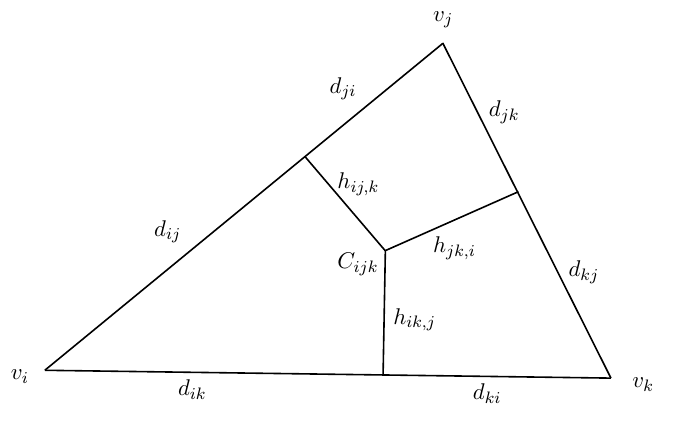}
  \caption{Signed distances of geometric center}
  \label{Signed distance of geometric center}
\end{figure}
Here the power distance of a point $p$ to the vertex $i$ is defined to be $\pi_p(i)=d^2(i,p)-\varepsilon_ir_i^2$,
where $d(i,p)$ is the Euclidean distance between $p$ and the vertex $i$. Please refer to Figure \ref{Signed distance of geometric center}
for the geometric center.
Denote $h_{jk,i}$ as the signed distance of the geometric center $C_{ijk}$ to the edge $\{jk\}$,
which is defined to be positive if $C_{ijk}$ is
on the same side of the line determined by $\{jk\}$ as the triangle $\{ijk\}$
and negative otherwise (or zero if $C_{ijk}$ is on the line).
Projections of $C_{ijk}$ to the edges $\{ij\}, \{ik\},\{jk\}$
give rise to the geometric centers of these edges, which are denoted by $C_{ij}, C_{ik}, C_{jk}$ respectively.
The signed distance $d_{ij}$ of $C_{ij}$ to the vertex $i$ is defined to be positive if $C_{ij}$ is on the same side as
$j$ along the line determined by $\{ij\}$ and negative otherwise (or zero if $C_{ij}$ is the same as $i$).
The signed distance $d_{ji}$ is defined similarly. Note that $d_{ij}+d_{ji}=l_{ij}$ and $d_{ij}\neq d_{ji}$ in general.
For nondegenerate radius vectors, we have
$
h_{jk,i}=\frac{d_{ji}-d_{jk}\cos \theta_j}{\sin \theta_j},\ \ d_{ij}=\frac{\varepsilon_ir_i^2+\eta_{ij}r_ir_j}{l_{ij}},
$
where $\theta_j$ is the inner angle at the vertex $j$ of the triangle $\{ijk\}$.
By direct calculations, we have
\begin{equation}\label{relationship of h_jki and h_i}
\begin{aligned}
h_{jk,i}=\frac{r_i^2r_j^2r_k^2}{Al_{jk}}\kappa_ih_i,
\end{aligned}
\end{equation}
where $A=l_{ij}l_{ik}\sin \theta_{i}$.
Lemma \ref{sign of hi hj hk} implies that the geometric center $C_{ijk}$  does not lie in some
region in the plane determined by the triangle as the nondegenerate radius vector tends to be degenerate.
Note that $h_i,h_j,h_k$ are defined for all radius vectors in $\mathbb{R}^3_{>0}$,
while $h_{ij,k}, h_{ik,j}, h_{jk,i}$ are defined only for nondegenerate radius vectors.
\end{remark}

Now we can give an analytic characterization of the admissible space of nondegenerate Euclidean discrete conformal factors on $(\{ijk\}, \varepsilon,\eta)$.
The main result is as follows.

\begin{theorem}\label{theorem admissible space Euclidean}
For the weighted triangle $(\{ijk\}, \varepsilon,\eta)$, the admissible space $\Omega_{ijk}^E(\eta)$ of
nondegenerate radius vectors
is a nonempty simply connected open set whose boundary components are analytic.
Furthermore,
$$
\Omega_{ijk}^E(\eta)=\mathbb{R}^3_{>0}\setminus\sqcup_{\alpha\in \Lambda}V_\alpha,
$$
where
$\Lambda=\{q\in \{i,j,k\}|A_q=\eta_{st}^2-\varepsilon_s\varepsilon_t>0, \{q, s, t\}=\{i, j, k\}\}$,
$\sqcup_{\alpha\in \Lambda}V_\alpha$ is a disjoint union of $V_\alpha$ and
$V_\alpha$ is a closed region in $\mathbb{R}^3_{>0}$ bounded by an analytic function defined on $\mathbb{R}^2_{>0}$ by
\begin{equation}\label{discription of V_i}
\begin{aligned}
V_i=\{(r_i, r_j, r_k)\in \mathbb{R}^3_{>0}| \kappa_i\geq \frac{-B_i+\sqrt{\Delta_i}}{2A_i}\}
   =\{(r_i, r_j, r_k)\in \mathbb{R}^3_{>0}| r_i\leq \frac{2A_i}{-B_i+\sqrt{\Delta_i}}\}.
\end{aligned}
\end{equation}
\end{theorem}

To prove Theorem \ref{theorem admissible space Euclidean}, we first prove the following result.

\begin{lemma}\label{A_i nonpositive}
For the weighted triangle $(\{ijk\}, \varepsilon,\eta)$,
if $\varepsilon_j\varepsilon_k-\eta_{jk}^2\geq 0$, $\varepsilon_i\varepsilon_k-\eta_{ik}^2\geq 0$ and $\varepsilon_i\varepsilon_j-\eta_{ij}^2\geq 0$,
then the admissible space $\Omega^E_{ijk}(\eta)$  in the parameter $r$ is $\mathbb{R}^3_{>0}$
and hence simply connected.
\end{lemma}

\proof
By Lemma \ref{lemma triangle equality euqi to Q>0}, we just need to prove that for any $r=(r_i, r_j, r_k)\in \mathbb{R}^3_{>0}$,
we have $Q^E>0$.
If $\varepsilon_j\varepsilon_k-\eta_{jk}^2\geq 0$, $\varepsilon_i\varepsilon_k-\eta_{ik}^2\geq 0$ and $\varepsilon_i\varepsilon_j-\eta_{ij}^2\geq 0$,
then we have $Q^E\geq 0$ by the definition (\ref{Q}) of $Q^E$ and the structure condition (\ref{structure condition 2}).
If $Q^E=0$, then $\eta_{jk}^2=\varepsilon_j\varepsilon_k, \eta_{ik}^2=\varepsilon_i\varepsilon_k, \eta_{ij}^2=\varepsilon_i\varepsilon_j$.
Combining this with the structure condition (\ref{structure condition 1}), we have $\eta_{jk}=\varepsilon_j\varepsilon_k, \eta_{ik}=\varepsilon_i\varepsilon_k, \eta_{ij}=\varepsilon_i\varepsilon_j$.
By the structure condition (\ref{structure condition 1}) again, this implies
$\varepsilon_i=\varepsilon_j=\varepsilon_k=\eta_{ij}=\eta_{ik}=\eta_{jk}=1$.
As a result, we have $Q^E=4\kappa_i\kappa_j+4\kappa_i\kappa_k+4\kappa_j\kappa_k>0$
for any $r=(r_i, r_j, r_k)\in \mathbb{R}^3_{>0}$. It is a contradiction.
Therefore, the admissible space $\Omega^E_{ijk}(\eta)=\mathbb{R}^3_{>0}$ and hence simply connected.
\qed

By Lemma \ref{A_i nonpositive}, we just need to study the case that at least one of $\varepsilon_j\varepsilon_k-\eta_{jk}^2$, $\varepsilon_i\varepsilon_k-\eta_{ik}^2$, $\varepsilon_i\varepsilon_j-\eta_{ij}^2$ is negative.
Suppose $r=(r_i, r_j, r_k)\in \mathbb{R}^3_{>0}$ is a degenerate radius vector.
Then $Q^E\leq 0$ by Lemma \ref{lemma triangle equality euqi to Q>0}.
By Lemma \ref{sign of hi hj hk}, one of $h_i, h_j, h_k$ is negative and the other two are positive.
Without loss of generality, assume $h_i<0$, $h_j>0, h_k>0$ at $r$.
By the definition (\ref{hi hj hk}) of $h_i$ and the structure condition (\ref{structure condition 2}), we have
$(\eta_{jk}^2-\varepsilon_j\varepsilon_k)\kappa_i>\gamma_k\kappa_j+\gamma_j\kappa_k\geq 0$.
This implies $\eta_{jk}^2-\varepsilon_j\varepsilon_k>0$.
Taking $Q^E$ as a quadratic function of $\kappa_i, \kappa_j, \kappa_k$. Then $Q^E\leq 0$ is equivalent to
\begin{equation}\label{Q in kappa i as quadratic}
\begin{aligned}
A_i\kappa_i^2+B_i\kappa_i+C_i\geq 0,
\end{aligned}
\end{equation}
where
\begin{equation}\label{A_i B_i C_i}
\begin{aligned}
A_i=&\eta_{jk}^2-\varepsilon_j\varepsilon_k>0,\\
B_i=&-2(\gamma_k\kappa_j+\gamma_j\kappa_k)\leq 0,\\
C_i=&(\eta_{ik}^2-\varepsilon_i\varepsilon_k)\kappa_j^2+(\eta_{ij}^2-\varepsilon_i\varepsilon_j)\kappa_k^{2}-2\kappa_j\kappa_k\gamma_i.
\end{aligned}
\end{equation}

\begin{lemma}\label{discriminant of kappa i}
For the weighted triangle $(\{ijk\}, \varepsilon,\eta)$,
if $A_i=\eta_{jk}^2-\varepsilon_j\varepsilon_k>0$,
then the discriminant $\Delta_i=B_i^2-4A_iC_i$ for (\ref{Q in kappa i as quadratic}) is positive.
\end{lemma}
\proof
By direct calculations, we have
\begin{equation}\label{Delta_i expression}
\begin{aligned}
\Delta_i
=4(\varepsilon_j\kappa_k^2+\varepsilon_k\kappa_j^2+2\eta_{jk}\kappa_j\kappa_k)
(\varepsilon_i\eta_{jk}^2+\varepsilon_j\eta_{ik}^2+\varepsilon_k\eta_{ij}^2+2\eta_{ij}\eta_{ik}\eta_{jk}-\varepsilon_i\varepsilon_j\varepsilon_k).
\end{aligned}
\end{equation}
By the structure condition (\ref{structure condition 1}) and the Cauchy inequality, we have
$\varepsilon_j\kappa_k^2+\varepsilon_k\kappa_j^2+2\eta_{jk}\kappa_j\kappa_k\geq
2(\varepsilon_j\varepsilon_k+\eta_{jk})\kappa_j\kappa_k>0$.
Therefore, the sign of $\Delta_i$ is determined by the term
$\varepsilon_i\eta_{jk}^2+\varepsilon_j\eta_{ik}^2+\varepsilon_k\eta_{ij}^2+2\eta_{ij}\eta_{ik}\eta_{jk}-\varepsilon_i\varepsilon_j\varepsilon_k$
in (\ref{Delta_i expression}), which is symmetric in $i,j,k$.

If one of $\varepsilon_i, \varepsilon_j, \varepsilon_k$ is zero, say $\varepsilon_i=0$, we have $\eta_{ij}>0$ and $\eta_{ik}>0$ by
the structure condition (\ref{structure condition 1}). This implies
$\varepsilon_i\eta_{jk}^2+\varepsilon_j\eta_{ik}^2+\varepsilon_k\eta_{ij}^2+2\eta_{ij}\eta_{ik}\eta_{jk}-\varepsilon_i\varepsilon_j\varepsilon_k
=(\varepsilon_j\eta_{ik}-\varepsilon_k\eta_{ij})^2+2\eta_{ij}\eta_{ik}(\eta_{jk}+\varepsilon_j\varepsilon_k)>0.$
Therefore, $\Delta_i>0$.

If $\varepsilon_i=\varepsilon_j=\varepsilon_k=1$, we have
$$\varepsilon_i\eta_{jk}^2+\varepsilon_j\eta_{ik}^2+\varepsilon_k\eta_{ij}^2+2\eta_{ij}\eta_{ik}\eta_{jk}-\varepsilon_i\varepsilon_j\varepsilon_k
=\eta_{jk}^2+\eta_{ik}^2+\eta_{ij}^2+2\eta_{ij}\eta_{ik}\eta_{jk}-1.$$
The positivity of this term under the condition $A_i=\eta_{jk}^2-1>0$ has been proved in Lemma 2.3 in \cite{X3}.
For completeness, we present a proof here.
If $\eta_{ij}\geq 0$ and $\eta_{ik}\geq 0$,
by the Cauchy inequality and the structure condition (\ref{structure condition 1}), we have
$\eta_{jk}^2+\eta_{ik}^2+\eta_{ij}^2+2\eta_{ij}\eta_{ik}\eta_{jk}-1\geq 2\eta_{ij}\eta_{ik}(\eta_{jk}+1)+\eta_{jk}^2-1\geq \eta_{jk}^2-1>0$.
If $\eta_{ij}<0$, then $\eta_{ij}\in (-1, 0)$ by the structure condition (\ref{structure condition 1}). This implies
$\eta_{jk}^2+\eta_{ik}^2+\eta_{ij}^2+2\eta_{ij}\eta_{ik}\eta_{jk}-1=(\eta_{ik}+\eta_{ij}\eta_{jk})^2+(\eta_{jk}^2-1)(1-\eta_{ij}^2)\geq (\eta_{jk}^2-1)(1-\eta_{ij}^2)>0$.
The same arguments apply to $\eta_{ik}<0$.
Therefore, $\Delta_i>0$.
 \qed

By the proof of Lemma \ref{discriminant of kappa i}, we have the following corollary.

\begin{corollary}\label{G positive}
For the weighted triangle $(\{ijk\}, \varepsilon,\eta)$,
if one of $\eta_{jk}^2-\varepsilon_j\varepsilon_k, \eta_{ik}^2-\varepsilon_i\varepsilon_k,\eta_{ij}^2-\varepsilon_i\varepsilon_j$
is positive, then the term
\begin{equation}\label{term G}
  G:=\varepsilon_i\eta_{jk}^2+\varepsilon_j\eta_{ik}^2+\varepsilon_k\eta_{ij}^2+2\eta_{ij}\eta_{ik}\eta_{jk}-\varepsilon_i\varepsilon_j\varepsilon_k
\end{equation}
is positive.
\end{corollary}

\begin{remark}\label{discriminant of kappa j k}
One can also take $Q^E$ as a quadratic function of $\kappa_j$ or $\kappa_k$ and define $\Delta_j$, $\Delta_k$
similarly.
By symmetry, we have
$\Delta_j>0$ if $\eta_{ik}^2-\varepsilon_i\varepsilon_k>0$ and
$\Delta_k>0$ if $\eta_{ij}^2-\varepsilon_i\varepsilon_j>0$.
\end{remark}

\textbf{Proof for Theorem \ref{theorem admissible space Euclidean}:}
We solve the admissible space of nondegenerate radius vectors for $(\{ijk\}, \varepsilon,\eta)$
by giving a precise description of the space of degenerate radius vectors.

Suppose $(r_i, r_j, r_k)\in \mathbb{R}^3_{>0}$ is a degenerate radius vector for $(\{ijk\}, \varepsilon,\eta)$.
By Lemma \ref{lemma triangle equality euqi to Q>0}, we have $Q^E=\kappa_ih_i+\kappa_jh_j+\kappa_kh_k\leq 0$.
By Lemma \ref{sign of hi hj hk}, one of $h_i,h_j,h_k$ is negative and the other two are positive.
Without loss of generality, assume $h_i<0, h_j>0, h_k>0$.
By $h_i<0$ and the structure condition (\ref{structure condition 2}), we have $A_i=\eta_{jk}^2-\varepsilon_j\varepsilon_k>0$.
Taking $Q^E\leq 0$ as a quadratic inequality of $\kappa_i$.
By Lemma \ref{discriminant of kappa i},  the solution of $Q^E\leq 0$, i.e. $A_i\kappa_i^2+B_i\kappa_i+C_i\geq 0$, is
\begin{equation*}
\begin{aligned}
\kappa_i\geq \frac{-B_i+\sqrt{\Delta_i}}{2A_i}\ \text{or} \ \kappa_i\leq \frac{-B_i-\sqrt{\Delta_i}}{2A_i}.
\end{aligned}
\end{equation*}
Note that
\begin{equation}\label{relation of h_i and Ai Bi}
\begin{aligned}
2A_i\kappa_i+B_i=2(\eta_{jk}^2-\varepsilon_j\varepsilon_k)\kappa_i-2(\gamma_k\kappa_j+\gamma_j\kappa_k)=-2h_i,
\end{aligned}
\end{equation}
we have $\kappa_i>\frac{-B_i}{2A_i}$ by $h_i<0$ and $A_i>0$.
This implies the solution $(r_i, r_j, r_k)\in \mathbb{R}^3_{>0}$ of $Q^E\leq 0$ with $h_i<0, h_j>0, h_k>0$ is $\kappa_i\geq \frac{-B_i+\sqrt{\Delta_i}}{2A_i}$.
Therefore, $\mathbb{R}^3_{>0}\setminus\Omega_{ijk}^E(\eta)\subset \cup_{\alpha\in \Lambda}V_\alpha$,
where
$\Lambda=\{q\in \{i,j,k\}|A_q=\eta_{st}^2-\varepsilon_s\varepsilon_t>0, \{q, s, t\}=\{i, j, k\}\}$, $V_i$ is defined by (\ref{discription of V_i})
and $V_j, V_k$ are defined similarly.

Conversely, suppose $(r_i, r_j, r_k)\in \cup_{\alpha\in \Lambda}V_\alpha\subseteq \mathbb{R}^3_{>0}$.
Without loss of generality, assume $(r_i, r_j, r_k)\in V_i$ and $A_i=\eta_{jk}^2-\varepsilon_j\varepsilon_k>0$.
Then $\kappa_i\geq \frac{-B_i+\sqrt{\Delta_i}}{2A_i}$ by the definition of $V_i$ in (\ref{discription of V_i}).
This is equivalent to $2A_i\kappa_i+B_i\geq \sqrt{\Delta_i}$
by $A_i>0$. Taking the square of both sides of this inequality gives $A_i\kappa_i^2+B_i\kappa_i+C_i\geq 0$, which is equivalent to $Q^E\leq 0$.
Therefore,
$\cup_{\alpha\in \Lambda}V_\alpha\subset\mathbb{R}^3_{>0}\setminus\Omega_{ijk}^E(\eta).$
In summary, we have
$\cup_{\alpha\in \Lambda}V_\alpha=\mathbb{R}^3_{>0}\setminus\Omega_{ijk}^E(\eta)$.

To see that $V_\alpha\cap V_\beta=\emptyset$ for distinct $\alpha$ and $\beta$ in $\Lambda$, suppose otherwise there exists some $(r_i, r_j, r_k)\in \mathbb{R}^3_{>0}$ with
$(r_i, r_j, r_k)\in V_i\cap V_j$. Then $A_i>0, A_j>0$.
By Lemma \ref{discriminant of kappa i} and Remark \ref{discriminant of kappa j k}, this implies
$\Delta_i>0, \Delta_j> 0$.
By $(r_i, r_j, r_k)\in V_i$, we have
$\kappa_i\geq \frac{-B_i+\sqrt{\Delta_i}}{2A_i}$.
Then by $A_i>0$ and (\ref{relation of h_i and Ai Bi}), we have
$h_i=-\frac{1}{2}(2A_i\kappa_i+B_i)\leq -\frac{1}{2}\sqrt{\Delta_i}<0$.
Lemma \ref{sign of hi hj hk} further implies $h_i<0, h_j>0, h_k>0$.
The same arguments applying to $(r_i, r_j, r_k)\in V_j$ shows that $h_j<0, h_i>0, h_k>0$. This is a contradiction. Therefore,
$V_\alpha\cap V_\beta=\emptyset$ for $\forall \alpha, \beta\in \Lambda, \alpha\neq\beta$.

Therefore,
$\Omega_{ijk}^E(\eta)=\mathbb{R}^3_{>0}\setminus\sqcup_{\alpha\in \Lambda}V_\alpha$.
As a result, the admissible space $\Omega_{ijk}^E(\eta)$ is homotopy equivalent to $\mathbb{R}^3_{>0}$ and hence simply connected.
\qed

\begin{remark}\label{sign of hi hj hj in Vi Euclidean}
By the proof of Theorem \ref{theorem admissible space Euclidean},
if $V_i$ defined by (\ref{discription of V_i}) is nonempty and $(r_i,r_j,r_k)\in V_i$, we have
$h_i<0, h_j>0,h_k>0$ at $(r_i,r_j,r_k)$.
\end{remark}

\begin{remark}
The method of characterizing the admissible space of nondegenerate discrete conformal factors on a weighted triangle in the proof of Theorem \ref{theorem admissible space Euclidean}
provides a unified approach to similar problems for other types of discrete conformal structures. See \cite{HX, X3, XZ1} for example.
The analytical characterization of the admissible space of nondegenerate discrete conformal factors on a weighted triangle
has some other applications.
See \cite{CLXZ} for example for some applications in the rigidity of infinite inversive distance circle packings on the plane and the convergence of the inversive distance circle packings.
\end{remark}

Define
\begin{equation*}
\begin{aligned}
\Omega_{ijk}^E=\{(r_i, r_j, r_k, \eta_{ij}, \eta_{ik}, \eta_{jk})\in \mathbb{R}^3_{>0}\times \mathbb{R}^3|
\eta \ \text{satisfies}\ (\ref{structure condition 1}), (\ref{structure condition 2}) \ \text{and } (r_i, r_j, r_k)\in \Omega_{ijk}^E(\eta)\}.
\end{aligned}
\end{equation*}
We call $\Omega_{ijk}^E$ as the parameterized admissible space of nondegenerate radius vectors for the triangle $\{ijk\}$.
The parameterized admissible space $\Omega_{ijk}^E$ contains some points with good properties.
\begin{lemma}\label{hi sign at special point}
The point $(r_i, r_j, r_k, \eta_{ij}, \eta_{ik}, \eta_{jk})=(1,1,1,1,1,1)$ is contained in $\Omega_{ijk}^E$.
Furthermore, $h_i>0, h_j>0, h_k>0$ at this point.
\end{lemma}
\proof
As $\varepsilon_i,\varepsilon_j,\varepsilon_k\in \{0,1\}$, it is straight forward to check that $(\eta_{ij}, \eta_{ik}, \eta_{jk})=(1,1,1)$
satisfies the structure conditions (\ref{structure condition 1}) and (\ref{structure condition 2}).
By $\varepsilon_i,\varepsilon_j,\varepsilon_k\in \{0,1\}$ and the definition (\ref{hi hj hk}) of $h_i, h_j, h_k$,  we have
\begin{equation*}
\begin{aligned}
h_i=&(\varepsilon_j\varepsilon_k-\eta_{jk}^2)\kappa_i+\kappa_j\gamma_k+\kappa_k\gamma_j=\varepsilon_j\varepsilon_k+\varepsilon_j+\varepsilon_k+1>0,\\
h_j=&(\varepsilon_i\varepsilon_k-\eta_{ik}^2)\kappa_j+\kappa_i\gamma_k+\kappa_k\gamma_i=\varepsilon_i\varepsilon_k+\varepsilon_i+\varepsilon_k+1>0,\\
h_k=&(\varepsilon_i\varepsilon_j-\eta_{ij}^2)\kappa_k+\kappa_i\gamma_j+\kappa_j\gamma_i=\varepsilon_i\varepsilon_j+\varepsilon_i+\varepsilon_j+1>0
\end{aligned}
\end{equation*}
at $(r_i, r_j, r_k, \eta_{ij}, \eta_{ik}, \eta_{jk})=(1,1,1,1,1,1)$,
which implies $Q^E=\kappa_ih_i+\kappa_jh_j+\kappa_kh_k>0$.
Therefore, by Lemma \ref{lemma triangle equality euqi to Q>0}, $(1,1,1,1,1,1)\in \Omega_{ijk}^E$. \qed

Theorem \ref{theorem admissible space Euclidean} have the following corollary on the parameterized admissible space $\Omega_{ijk}^E$.
\begin{corollary}\label{coro connectivity of para admi space}
For the triangle $\sigma=\{ijk\}$ with a weight $\varepsilon:V_\sigma\rightarrow \{0,1\}$,
 the parameterized admissible space $\Omega_{ijk}^E$ is connected.
\end{corollary}
\proof
Set
\begin{equation*}
\begin{aligned}
\Gamma=\{(\eta_{ij}, \eta_{ik}, \eta_{jk})\in \mathbb{R}^3|(\eta_{ij}, \eta_{ik}, \eta_{jk}) \ \text{satisfies}\  (\ref{structure condition 1}), (\ref{structure condition 2}) \}.
\end{aligned}
\end{equation*}
Then $\Omega_{ijk}^E$ is a fiber bundle over $\Gamma$,  and the fiber over $\eta=(\eta_{ij}, \eta_{ik}, \eta_{jk})\in \Gamma$ is the connected admissible space $\Omega_{ijk}^E(\eta)$.
We will prove that $\Gamma$ is path connected. As a result,  the connectivity of $\Omega_{ijk}^E$ follows
by Theorem \ref{theorem admissible space Euclidean} and
the continuity of $Q$ as a function of $(r_i, r_j, r_k, \eta_{ij}, \eta_{ik}, \eta_{jk})$.

It is obviously that $\mathbb{R}^3_{>0}\subset \Gamma$, which is path connected.
We will show that any point in $\Gamma$ can be connected to  $\mathbb{R}^3_{>0}$ by a path in $\Gamma$.
As the boundary of $\mathbb{R}^3_{>0}$ is connected to  $\mathbb{R}^3_{>0}$, we just need to consider the case that some
component of $(\eta_{ij}, \eta_{ik}, \eta_{jk})\in \Gamma$ is negative.
Without loss of generality, assume $\eta_{ij}<0$, then $\varepsilon_i=\varepsilon_j=1$ by the structure condition $\eta_{ij}+\varepsilon_i\varepsilon_j>0$. Therefore, we just need to consider the cases $\varepsilon_k=0$ and $\varepsilon_k=1$.


In the case of $\varepsilon_k=0$, the structure conditions (\ref{structure condition 1}), (\ref{structure condition 2}) are equivalent to
\begin{equation}\label{structure condition 1 proof a}
\begin{aligned}
1+\eta_{ij}>0, \eta_{ik}>0, \eta_{jk}>0
\end{aligned}
\end{equation}
and
\begin{equation}\label{structure condition 2 proof a}
\begin{aligned}
\eta_{jk}+\eta_{ij}\eta_{ik}\geq0, \eta_{ik}+\eta_{ij}\eta_{jk}\geq0, \eta_{ik}\eta_{jk}\geq0.
\end{aligned}
\end{equation}
If $(\eta_{ij}, \eta_{ik}, \eta_{jk})\in \Gamma$ and $\eta_{ij}<0$,
it is straightforward to check that $(t\eta_{ij}, \eta_{ik}, \eta_{jk})$ satisfies
(\ref{structure condition 1 proof a}) and (\ref{structure condition 2 proof a}) for any $t\in [0, 1]$.
This implies $(t\eta_{ij}, \eta_{ik}, \eta_{jk})\in \Gamma$, $\forall t\in [0, 1]$,
which is a path connecting $(\eta_{ij}, \eta_{ik}, \eta_{jk})$ and $\mathbb{R}^3_{\geq 0}$.
Therefore, $\Gamma$ is path connected.

In the case of $\varepsilon_k=1$, the structure conditions (\ref{structure condition 1}), (\ref{structure condition 2}) are equivalent to
\begin{equation}\label{structure condition 1 proof b}
\begin{aligned}
1+\eta_{ij}>0, 1+\eta_{ik}>0, 1+\eta_{jk}>0
\end{aligned}
\end{equation}
and
\begin{equation}\label{structure condition 2 proof b}
\begin{aligned}
\eta_{jk}+\eta_{ij}\eta_{ik}\geq0, \eta_{ik}+\eta_{ij}\eta_{jk}\geq0, \eta_{ij}+\eta_{ik}\eta_{jk}\geq0.
\end{aligned}
\end{equation}
In this case, the path connectivity of $\Gamma$ has been proved in \cite{X3}. For completeness, we present the proof here.
By the structure conditions (\ref{structure condition 1 proof b}) and taking the sum of the equations in (\ref{structure condition 2 proof b}) in pairs, we have $\eta_{ij}+\eta_{ik}\geq 0$,
$\eta_{ij}+\eta_{jk}\geq 0$, $\eta_{ik}+\eta_{jk}\geq 0$. This implies at most one of $\eta_{ij}, \eta_{ik}, \eta_{jk}$
is negative. By the assumption that $\eta_{ij}<0$, we have $\eta_{ik}> 0, \eta_{jk}> 0$.
It is straightforward to check that $(t\eta_{ij}, \eta_{ik}, \eta_{jk})$
satisfies
(\ref{structure condition 1 proof b}) and (\ref{structure condition 2 proof b}) for any $t\in [0, 1]$.
This implies $(t\eta_{ij}, \eta_{ik}, \eta_{jk})\in \Gamma$, $\forall t\in [0, 1]$,
which is a path connecting $(\eta_{ij}, \eta_{ik}, \eta_{jk})$ and $\mathbb{R}^3_{\geq 0}$.
Therefore, $\Gamma$ is path connected.
\qed

\subsection{Negative semi-definiteness of the Jacobian matrix in the Euclidean background geometry}

Let $(\{ijk\}, \varepsilon,\eta)$ be a nondegenerate weighted Euclidean triangle with edge lengths given by (\ref{defn of Euclidean length}).
And $\theta_i, \theta_j, \theta_k$ are the inner angles at the vertices $i, j, k$ respectively.
Set $u_i=f_i=\ln r_i$.
\begin{lemma}[\cite{G3}]\label{Euclidean symmetry}
Let $(\{ijk\}, \varepsilon,\eta)$ be a weighted triangle and $(r_i, r_j, r_k)\in \mathbb{R}^3_{>0}$ is a nondegenerate radius vector on $(\{ijk\}, \varepsilon,\eta)$.
Then
\begin{equation}\label{derivative of theta i of uj Euclidean}
\begin{aligned}
\frac{\partial \theta_i}{\partial u_j}=\frac{\partial \theta_j}{\partial u_i}
=\frac{r_i^2r_j^2r_k^2}{Al_{ij}^2}[(\varepsilon_i\varepsilon_j-\eta_{ij}^2)\kappa_k^2+\gamma_i\kappa_j\kappa_k+\gamma_j\kappa_i\kappa_k]
=\frac{r_i^2r_j^2r_k}{Al_{ij}^2}h_k
\end{aligned}
\end{equation}
and
\begin{equation}\label{derivative of theta i of ui Euclidean}
\begin{aligned}
\frac{\partial \theta_i}{\partial u_i}
=-\frac{\partial \theta_i}{\partial u_j}-\frac{\partial \theta_i}{\partial u_k},
\end{aligned}
\end{equation}
where $A=l_{ij}l_{ik}\sin \theta_i$.
\end{lemma}
\proof
By the chain rules, we have
\begin{equation}\label{derivative of theta wrt uj}
\begin{aligned}
\frac{\partial \theta_i}{\partial u_j}
=&\frac{\partial \theta_i}{\partial l_{jk}}\frac{\partial l_{jk}}{\partial u_j}
  +\frac{\partial \theta_i}{\partial l_{ik}}\frac{\partial l_{ik}}{\partial u_j}
  +\frac{\partial \theta_i}{\partial l_{ij}}\frac{\partial l_{ij}}{\partial u_j}.
\end{aligned}
\end{equation}
By the derivative cosine law (\cite{CL}, Lemma A1), we have
\begin{equation}\label{Euclidean derivative cosine law}
\begin{aligned}
\frac{\partial \theta_i}{\partial l_{jk}}=\frac{l_{jk}}{A},
\frac{\partial \theta_i}{\partial l_{ik}}=\frac{-l_{jk}\cos\theta_k}{A},
\frac{\partial \theta_i}{\partial l_{ij}}=\frac{-l_{jk}\cos\theta_j}{A},
\end{aligned}
\end{equation}
where $A=l_{ik}l_{jk}\sin \theta_k$.
By the definition (\ref{definition of Euclidean length}) of $l_{ij}, l_{ik}, l_{jk}$ in $r_i, r_j, r_k$, we have
\begin{equation}\label{djk dji proof}
\begin{aligned}
\frac{\partial l_{jk}}{\partial u_j}
=\frac{\varepsilon_jr_j^2+\eta_{jk}r_jr_k}{l_{jk}},\ \
\frac{\partial l_{ik}}{\partial u_j}=0,\ \
\frac{\partial l_{ij}}{\partial u_j}
=\frac{\varepsilon_jr_j^2+\eta_{ij}r_ir_j}{l_{ij}}.
\end{aligned}
\end{equation}
Submitting (\ref{Euclidean derivative cosine law}) and (\ref{djk dji proof}) into (\ref{derivative of theta wrt uj}),
we have
\begin{equation}\label{derivative of theta i of uj}
\begin{aligned}
\frac{\partial \theta_i}{\partial u_j}
=&\frac{l_{jk}}{A}\cdot\frac{\varepsilon_jr_j^2+\eta_{jk}r_jr_k}{l_{jk}}
  +\frac{-l_{jk}\cos\theta_j}{A}\cdot\frac{\varepsilon_jr_j^2+\eta_{ij}r_ir_j}{l_{ij}}\\
=&\frac{1}{2Al_{ij}^2}[2(\varepsilon_jr_j^2+\eta_{jk}r_jr_k)l_{ij}^2+(l_{ik}^2-l_{ij}^2-l_{jk}^2)(\varepsilon_jr_j^2+\eta_{ij}r_ir_j)]\\
=&\frac{r_i^2r_j^2r_k^2}{Al_{ij}^2}[(\varepsilon_i\varepsilon_j-\eta_{ij}^2)\kappa_k^2+\gamma_i\kappa_j\kappa_k+\gamma_j\kappa_i\kappa_k]\\
=&\frac{r_i^2r_j^2r_k}{Al_{ij}^2}h_k,
\end{aligned}
\end{equation}
where the cosine law is used in the second line and the definition (\ref{definition of Euclidean length})
of edge lengths is used in the third line.
As the last line of (\ref{derivative of theta i of uj})
is symmetric in $i$ and $j$, we have
$\frac{\partial \theta_i}{\partial u_j}=\frac{\partial \theta_j}{\partial u_i}$.
Similarly, we have $\frac{\partial \theta_i}{\partial u_k}=\frac{\partial \theta_k}{\partial u_i}$.
The formula
$\frac{\partial \theta_i}{\partial u_i}
=-\frac{\partial \theta_i}{\partial u_j}-\frac{\partial \theta_i}{\partial u_k}$
follows from $\theta_i+\theta_j+\theta_k=\pi$, $\frac{\partial \theta_i}{\partial u_j}=\frac{\partial \theta_j}{\partial u_i}$ and
$\frac{\partial \theta_i}{\partial u_k}=\frac{\partial \theta_k}{\partial u_i}$.
\qed

\begin{remark}\label{geometric explaination of part theta part u}
The property $\frac{\partial \theta_i}{\partial u_j}=\frac{\partial \theta_j}{\partial u_i}$
in Lemma \ref{Euclidean symmetry} was proved by Glickenstein \cite{G3}.
Here we give a proof by direct calculations for completeness.
Combining (\ref{relationship of h_jki and h_i}) and Lemma \ref{Euclidean symmetry}, we have
\begin{equation}\label{geometric explain of deriv theta wrt u}
\begin{aligned}
\frac{\partial \theta_i}{\partial u_j}=\frac{h_{ij,k}}{l_{ij}}.
\end{aligned}
\end{equation}
The formula (\ref{geometric explain of deriv theta wrt u}) provides a nice geometric explanation for the derivative $\frac{\partial \theta_i}{\partial u_j}$, please refer to \cite{G3} for more information for this.
\end{remark}

\begin{remark}\label{derivative tends infty Euclidean}
By (\ref{derivative of theta i of uj Euclidean}), (\ref{derivative of theta i of ui Euclidean}) and Remark \ref{sign of hi hj hj in Vi Euclidean},
if  $(r_i, r_j, r_k)\in \Omega_{ijk}^E(\eta)$ tends to a point $(\overline{r}_i, \overline{r}_j, \overline{r}_k)\in \partial V_i$ with $V_i\neq \emptyset$,
 we have $\frac{\partial \theta_i}{\partial u_j}\rightarrow +\infty$, $\frac{\partial \theta_i}{\partial u_k}\rightarrow +\infty$ and
 $\frac{\partial \theta_i}{\partial u_i}\rightarrow -\infty$.
\end{remark}

Lemma \ref{Euclidean symmetry} shows that the Jacobian matrix
\begin{equation*}
\begin{aligned}
\Lambda^E_{ijk}:=\frac{\partial (\theta_i, \theta_j, \theta_k)}{\partial ( u_i, u_j, u_k)}
=\left(
   \begin{array}{ccc}
     \frac{\partial \theta_i}{\partial u_i} & \frac{\partial \theta_i}{\partial u_j} & \frac{\partial \theta_i}{\partial u_k} \\
     \frac{\partial \theta_j}{\partial u_i} & \frac{\partial \theta_j}{\partial u_j} & \frac{\partial \theta_j}{\partial u_k} \\
     \frac{\partial \theta_k}{\partial u_i} & \frac{\partial \theta_k}{\partial u_j} & \frac{\partial \theta_k}{\partial u_k} \\
   \end{array}
 \right)
\end{aligned}
\end{equation*}
is symmetric with $\{t(1,1,1)^T|t\in \mathbb{R}\}$ in its kernel.
Furthermore, We have the following result on the rank of the Jacobian matrix  $\Lambda^E_{ijk}$.

\begin{lemma}\label{Euclidean rank lemma}
For the weighted triangle $(\{ijk\}, \varepsilon,\eta)$, the rank of $\Lambda^E_{ijk}$ is $2$
for any nondegenerate radius vector.
\end{lemma}

\proof
By the chain rules, we have
\begin{equation}\label{Euclidean chain rules}
\begin{aligned}
\frac{\partial (\theta_i, \theta_j, \theta_k)}{\partial ( u_i, u_j, u_k)}
=\frac{\partial (\theta_i, \theta_j, \theta_k)}{\partial ( l_{jk}, l_{ik}, l_{ij})}
\cdot\frac{\partial ( l_{jk}, l_{ik}, l_{ij})}{\partial ( u_i, u_j, u_k)}.
\end{aligned}
\end{equation}
By the derivative cosine law (\cite{CL}, Lemma A1), we have
\begin{equation*}
\begin{aligned}
\frac{\partial (\theta_i, \theta_j, \theta_k)}{\partial ( l_{jk}, l_{ik}, l_{ij})}
=\frac{1}{A}\left(
              \begin{array}{ccc}
                l_{jk} &   &   \\
                  & l_{ik} &   \\
                  &   & l_{ij} \\
              \end{array}
            \right)
            \left(
              \begin{array}{ccc}
                1 & -\cos\theta_k & -\cos\theta_j \\
                -\cos\theta_k & 1 & -\cos\theta_i \\
                -\cos\theta_j & -\cos\theta_i & 1 \\
              \end{array}
            \right).
\end{aligned}
\end{equation*}
This matrix has rank $2$ and kernel $\{t(l_{jk}, l_{ik}, l_{ij})|t\in \mathbb{R}\}$
for $(l_{jk}, l_{ik}, l_{ij})$ satisfying the triangle inequalities.

Note that $d_{ij}=\frac{\partial l_{ij}}{\partial u_i}=\frac{\varepsilon_ir_i^2+\eta_{ij}r_ir_j}{l_{ij}}$.
By direct calculations,
\begin{equation*}
\begin{aligned}
\frac{\partial ( l_{jk}, l_{ik}, l_{ij})}{\partial ( u_i, u_j, u_k)}
=&\left(
   \begin{array}{ccc}
     0 & d_{jk} & d_{kj} \\
     d_{ik} & 0 & d_{ki} \\
     d_{ij} & d_{ji} & 0 \\
   \end{array}
 \right)\\
=&\left(
    \begin{array}{ccc}
      l_{jk}^{-1} &   &   \\
        & l_{ik}^{-1} &   \\
        &   & l_{ij}^{-1} \\
    \end{array}
  \right)\\
& \cdot
  \left(
    \begin{array}{ccc}
      0 & \varepsilon_jr_j+\eta_{jk}r_k & \varepsilon_kr_k+\eta_{jk}r_j \\
      \varepsilon_ir_i+\eta_{ik}r_k & 0 & \varepsilon_kr_k+\eta_{ik}r_i \\
      \varepsilon_ir_i+\eta_{ij}r_j & \varepsilon_jr_j+\eta_{ij}r_i & 0 \\
    \end{array}
  \right)
  \left(
    \begin{array}{ccc}
      r_i &   &   \\
        & r_j &   \\
        &   & r_k \\
    \end{array}
  \right).
\end{aligned}
\end{equation*}
This implies
\begin{equation}\label{Euclidean positive determinant}
\begin{aligned}
\det \frac{\partial ( l_{jk}, l_{ik}, l_{ij})}{\partial ( u_i, u_j, u_k)}
=&\frac{r_ir_jr_k}{l_{ij}l_{ik}l_{jk}}
  [2(\varepsilon_i\varepsilon_j\varepsilon_k+\eta_{ij}\eta_{ik}\eta_{jk}) r_ir_jr_k
     +r_i\gamma_i(\varepsilon_jr_j^2+\varepsilon_kr_k^2)\\
     &+r_j\gamma_j(\varepsilon_ir_i^2+\varepsilon_kr_k^2)
     +r_k\gamma_k(\varepsilon_ir_i^2+\varepsilon_jr_j^2)]\\
\geq &\frac{2r_i^2r_j^2r_k^2}{l_{ij}l_{ik}l_{jk}}[\varepsilon_i\varepsilon_j\varepsilon_k+\eta_{ij}\eta_{ik}\eta_{jk}
       +\gamma_i\varepsilon_j\varepsilon_k+\gamma_j\varepsilon_i\varepsilon_k+\gamma_k\varepsilon_i\varepsilon_j]\\
=&\frac{2r_i^2r_j^2r_k^2}{l_{ij}l_{ik}l_{jk}}(\varepsilon_i\varepsilon_j+\eta_{ij})
   (\varepsilon_i\varepsilon_k+\eta_{ik})(\varepsilon_j\varepsilon_k+\eta_{jk})\\
>&0,
\end{aligned}
\end{equation}
where the structure condition (\ref{structure condition 2}) is used in the second line and
the structure condition (\ref{structure condition 1}) is used in the last line.
The inequality (\ref{Euclidean positive determinant}) implies that $\frac{\partial ( l_{jk}, l_{ik}, l_{ij})}{\partial ( u_i, u_j, u_k)}$
is nonsingular.

By (\ref{Euclidean chain rules}), we have
the rank of $\Lambda^E_{ijk}=\frac{\partial (\theta_i, \theta_j, \theta_k)}{\partial ( u_i, u_j, u_k)}$ is $2$
for any nondegenerate radius vector on $(\{ijk\}, \varepsilon,\eta)$.
\qed


\begin{theorem}\label{Euclidean Jacobian negativity}
For the weighted triangle $(\{ijk\}, \varepsilon,\eta)$,
the Jacobian matrix $\Lambda^E_{ijk}=\frac{\partial (\theta_i, \theta_j, \theta_k)}{\partial ( u_i, u_j, u_k)}$
is negative semi-definite with rank $2$ and  has kernel $\{t(1,1,1)^T|t\in \mathbb{R}\}$ for any nondegenerate
Euclidean discrete conformal factor on $(\{ijk\}, \varepsilon,\eta)$.
\end{theorem}
\proof
By Lemma \ref{Euclidean rank lemma}, the matrix $\Lambda^E_{ijk}$
has two nonzero eigenvalues and one zero eigenvalue.
By the continuity of the eigenvalues of $\Lambda^E_{ijk}$
as functions of $(r_i, r_j, r_k, \eta_{ij}, \eta_{ik}, \eta_{jk})\in \Omega_{ijk}^E$ and
the connectivity of parameterized admissible space $\Omega_{ijk}^E$ in Corollary \ref{coro connectivity of para admi space}, to prove
$\Lambda^E_{ijk}$ is negative semi-definite,
we just need to prove $\Lambda^E_{ijk}$ is negative semi-definite with rank $2$
at some point in $\Omega_{ijk}^E$.
By Lemma \ref{hi sign at special point},
$h_i>0, h_j>0, h_k>0$ at the point
$(r_i, r_j, r_k, \eta_{ij}, \eta_{ik}, \eta_{jk})=(1,1,1,1,1,1)\in \Omega_{ijk}^E$.
By (\ref{relationship of h_jki and h_i}) and (\ref{geometric explain of deriv theta wrt u}),
this implies $\frac{\partial \theta_i}{\partial u_j}$,
$\frac{\partial \theta_i}{\partial u_k}$, $\frac{\partial \theta_j}{\partial u_k}$ are positive.
Then by the following well-known result from linear algebra,
$-\Lambda^E_{ijk}$ is positive semi-definite with rank $2$ and has kernel $\{t(1,1,1)^T|t\in \mathbb{R}\}$ at
$(r_i, r_j, r_k, \eta_{ij}, \eta_{ik}, \eta_{jk})=(1,1,1,1,1,1)\in \Omega_{ijk}^E$.
\begin{lemma}\label{diagonal dominant}
Suppose $A=[a_{ij}]_{n\times n}$ is a symmetric matrix.
\begin{description}
  \item[(a)] If $a_{ii}>\sum_{j\neq i}|a_{ij}|$ for all indices $i$, then $A$ is positive definite.
  \item[(b)] If $a_{ii}>0$ and $a_{ij}< 0$ for all $i\neq j$ so that $\sum_{i=1}^na_{ij}=0$ for all $j$,
  then $A$ is positive semi-definite so that its kernel is 1-dimensional.
\end{description}
\end{lemma}
One can refer to \cite{CL} for a proof of Lemma \ref{diagonal dominant}.
Therefore, $\Lambda^E_{ijk}$ is negative semi-definite  with rank $2$ and has kernel $\{t(1,1,1)^T|t\in \mathbb{R}\}$
for any point $(r_i, r_j, r_k, \eta_{ij}, \eta_{ik}, \eta_{jk})\in \Omega_{ijk}^E$. \qed

In the literature, the proof for the nonnegative semi-definiteness of the Jacobian matrix $\Lambda^E_{ijk}$
is based on direct and tedious calculations. See \cite{Guo,X1} for example.
The proof  of Theorem \ref{Euclidean Jacobian negativity}
based on parameterized admissible space provides a much simpler approach for such problems.

As a corollary of Theorem \ref{Euclidean Jacobian negativity}, we have the following result on the
Jacobian matrix $\Lambda^E=\frac{\partial (K_1,\cdots,K_N)}{\partial (u_1,\cdots, u_N)}$.

\begin{corollary}\label{Euclidean curvature Jacobian positivity}
Suppose $(M, \mathcal{T}, \varepsilon, \eta)$ is a weighted triangulated surface with
the weights $\varepsilon: V\rightarrow \{0, 1\}$ and $\eta: E\rightarrow \mathbb{R}$
satisfying the structure conditions (\ref{structure condition 1}) and (\ref{structure condition 2}).
Then the Jacobian matrix $\Lambda^E=\frac{\partial (K_1,\cdots,K_N)}{\partial (u_1,\cdots, u_N)}$
is symmetric and positive semi-definite with rank $N-1$ and has kernel $\{t\mathbf{1}\in \mathbb{R}^N|t\in \mathbb{R}\}$
for all nondegenerate Euclidean discrete conformal factors on $(M, \mathcal{T}, \varepsilon, \eta)$.
\end{corollary}
\proof
This follows from Theorem \ref{Euclidean Jacobian negativity} and the fact that $\Lambda^E=-\sum_{\{ijk\}\in F}\Lambda^E_{ijk}$,
where $\Lambda^E_{ijk}$ is extended by zeros to be an $N\times N$ matrix so that
$\Lambda^E_{ijk}$ acts on a vector $(v_1, \cdots, v_N)$ only on the coordinates
corresponding to vertices $v_i$, $v_j$ and $v_k$ in the triangle $\{ijk\}$.
\qed

\begin{remark}
Under an additional condition that the signed distance of geometric center to the edges are all positive for any triangle $\{ijk\}\in F$,
Glickenstein \cite{G4} and Glickenstein-Thomas \cite{GT} proved the positive semi-definiteness of
the Jacobian matrix $\Lambda^E=\frac{\partial (K_1,\cdots,K_N)}{\partial (u_1,\cdots, u_N)}$.
Corollary \ref{Euclidean curvature Jacobian positivity} generalizes Glickenstein-Thomas's result in that
it allows some of the signed distance to be negative.
For example, in the case that $\varepsilon\equiv 1$ and $\eta\equiv 2$,
if $r: V\rightarrow (0, +\infty)$  is a map with
$r\equiv 1$ except $r_i=1/5$ for some vertex $i\in V$,
then $r$ is a nondegenerate Euclidean radius vector on $(M, \mathcal{T}, \varepsilon, \eta)$.
By Corollary \ref{Euclidean curvature Jacobian positivity}, $\Lambda^E$ is positive semi-definite at $r$.
However, we have $h_i<0, h_j>0, h_k>0$ for any triangle $\{ijk\}$ at $i$, which implies
$h_{jk,i}<0, h_{ik,j}>0, h_{ij, k}>0$ at $r$ by (\ref{relationship of h_jki and h_i}).
\end{remark}

\subsection{Rigidity of Euclidean discrete conformal structures}
By Theorem \ref{theorem admissible space Euclidean} and Lemma \ref{Euclidean symmetry}, the following function
\begin{equation}\label{Euclidean Ricci energy function for triangle}
\begin{aligned}
\mathcal{E}_{ijk}(u_i,u_j,u_k)=\int_{(\overline{u}_i, \overline{u}_j, \overline{u}_k)}^{(u_i, u_j, u_k)}\theta_idu_i+\theta_jdu_j+\theta_kdu_k
\end{aligned}
\end{equation}
is a well-defined smooth function on $\Omega^E_{ijk}(\eta)$ with $\nabla_{u_i}\mathcal{E}_{ijk}=\theta_i$ and
$\mathcal{E}_{ijk}(u_i+t,u_j+t,u_k+t)=\mathcal{E}_{ijk}(u_i,u_j,u_k)+t\pi$.
The function $\mathcal{E}_{ijk}(u_i,u_j,u_k)$ is called the Ricci energy function for the weighted triangle $(\{ijk\}, \varepsilon,\eta)$.
It was first constructed by Glickenstein \cite{G3} for Glickenstein's Euclidean discrete conformal structures
under the assumption that the domain is simply connected.
Furthermore, Glickenstein-Thomas \cite{GT} used the Ricci energy function to prove a result on the local rigidity of Glickenstein's
Euclidean discrete conformal structures.
For completeness, we give a sketch of  Glickenstein-Thomas's arguments here.
By Theorem \ref{Euclidean Jacobian negativity}, $\mathcal{E}_{ijk}(u_i,u_j,u_k)$ is a locally concave function defined on $\Omega^E_{ijk}(\eta)$.
Set
\begin{equation}\label{Euclidean Ricci energy function}
\begin{aligned}
\mathcal{E}(u_1, \cdots, u_N)=2\pi\sum_{i\in V}u_i-\sum_{\{ijk\}\in F}\mathcal{E}_{ijk}(u_i,u_j,u_k).
\end{aligned}
\end{equation}
We call $\mathcal{E}(u_1, \cdots, u_N)$  as the Ricci energy function for $(M, \mathcal{T}, \varepsilon, \eta)$. It is defined on the admissible space $\Omega^E$ of nondegenerate Euclidean discrete conformal factors.
By Corollary \ref{Euclidean curvature Jacobian positivity}, $\mathcal{E}$ is a locally convex function defined on $\Omega^E$ with
$\mathcal{E}(u_1+t, \cdots, u_N+t)=\mathcal{E}(u_1, \cdots, u_N)+2t\pi\chi(M)$ and $\nabla_{u_i} \mathcal{E}=K_i$.
The local rigidity of Glickenstein's Euclidean discrete conformal structures follows by the following well-known result from analysis.

\begin{lemma}\label{injectivity of convex function}
If $W: \Omega\rightarrow \mathbb{R}$ is a $C^2$-smooth strictly convex function defined on a convex domain $\Omega\subseteq \mathbb{R}^n$,
then its gradient $\nabla W: \Omega\rightarrow \mathbb{R}^n$ is injective.
\end{lemma}
To prove the global rigidity of Glickenstein's Euclidean discrete conformal structures,
we need to extend the inner angles of a triangle $\{ijk\}$
defined for nondegenerate radius vectors
to be a globally defined function for all radius vectors $(r_i, r_j, r_k)\in \mathbb{R}^3_{>0}$.

\begin{lemma}\label{Euclidean extension}
For the weighted triangle $(\{ijk\}, \varepsilon,\eta)$,
the inner angles $\theta_i, \theta_j, \theta_k$
defined for nondegenerate radius vectors
can be extended by constants to be continuous functions $\widetilde{\theta}_i, \widetilde{\theta}_j, \widetilde{\theta}_k$
defined for $(r_i, r_j, r_k)\in \mathbb{R}^3_{>0}$ by setting
\begin{equation}\label{extension of theta_i}
\begin{aligned}
\widetilde{\theta}_i(r_i, r_j, r_k)=\left\{
                                      \begin{array}{ll}
                                        \theta_i, & \hbox{if $(r_i, r_j, r_k)\in \Omega^E_{ijk}(\eta)$;} \\
                                        \pi, & \hbox{if $(r_i, r_j, r_k)\in V_i$;} \\
                                        0, & \hbox{otherwis.}
                                      \end{array}
                                    \right.
\end{aligned}
\end{equation}
\end{lemma}
\proof
By Theorem \ref{theorem admissible space Euclidean},
$\Omega_{ijk}^E(\eta)=\mathbb{R}^3_{>0}\setminus\sqcup_{\alpha\in \Lambda}V_\alpha$,
where
$\Lambda=\{q\in \{i,j,k\}|A_q=\eta_{st}^2-\varepsilon_s\varepsilon_t>0, \{q, s, t\}=\{i, j, k\}\}$ and
$V_\alpha$ is a closed region in $\mathbb{R}^3_{>0}$ bounded by the analytical function in (\ref{discription of V_i})
defined on $\mathbb{R}^2_{>0}$.

If $\Lambda= \varnothing$, then $\Omega_{ijk}^E(\eta)=\mathbb{R}^3_{>0}$ and
$\theta_i, \theta_j, \theta_k$ is defined for all $(r_i, r_j, r_k)\in \mathbb{R}^3_{>0}$.

If $\Lambda\neq \varnothing$, let $V_i$ be a connected component of $\mathbb{R}^3_{>0}\setminus \Omega_{ijk}^E(\eta)$.
Suppose $(r_i, r_j, r_k)\in \Omega_{ijk}^E(\eta)$ tends to a point
$(\overline{r}_i, \overline{r}_j, \overline{r}_k)$ in the boundary $\partial V_i$ of $V_i$ in $\mathbb{R}^3_{>0}$.
By Heron's formula, we have
\begin{equation}\label{Heron's formula}
\begin{aligned}
4l_{ij}^2l_{ik}^2\sin^2\theta_i
=(l_{ij}+l_{ik}+l_{jk})(l_{ij}+l_{ik}-l_{jk})(l_{ij}-l_{ik}+l_{jk})(-l_{ij}+l_{ik}+l_{jk})\rightarrow 0.
\end{aligned}
\end{equation}
Note that for any $r_i, r_j>0$, by the structure condition (\ref{structure condition 1}) and Cauchy inequality, we have
$\varepsilon_ir_i^2+\varepsilon_jr_j^2+2\eta_{ij}r_ir_j\geq 2(\varepsilon_i\varepsilon_j+\eta_{ij})r_ir_j>0$.
This implies $l_{ij}, l_{ik}$ tend to positive numbers as $(r_i, r_j, r_k)\rightarrow (\overline{r}_i, \overline{r}_j, \overline{r}_k)$.
Combining this and (\ref{Heron's formula}), we have $\sin\theta_i$ tends to zero.
Therefore, $\theta_i$ tends to $0$ or $\pi$. Similarly, we have $\theta_j, \theta_k$ tends to $0$ or $\pi$.

By Remark \ref{sign of hi hj hj in Vi Euclidean},  we have
 $h_i<0$, $h_j>0$ and $h_k>0$ at $(\overline{r}_i, \overline{r}_j, \overline{r}_k)\in \partial V_i$.
By the continuity of $h_i, h_j, h_k$,
there exists some neighborhood $U$
of $(\overline{r}_i, \overline{r}_j, \overline{r}_k)$ in $\mathbb{R}^3_{>0}$
such that $h_i<0$, $h_j>0$, $h_k>0$ for $(r_i, r_j, r_k)\in \Omega^E_{ijk}(\eta)\cap U$.
Combining $h_k>0$, (\ref{relationship of h_jki and h_i}) and (\ref{geometric explain of deriv theta wrt u}), we have
$
\frac{\partial \theta_i}{\partial u_j}=\frac{r_i^2r_j^2r_k^2}{Al^2_{ij}}\kappa_kh_k>0
$
for $(r_i, r_j, r_k)\in \Omega^E_{ijk}(\eta)\cap U$.
Similarly, we have $\frac{\partial \theta_i}{\partial u_k}>0$ for $(r_i, r_j, r_k)\in \Omega^E_{ijk}(\eta)\cap U$.
By Lemma \ref{Euclidean symmetry}, we have
$\frac{\partial \theta_i}{\partial u_i}
=-\frac{\partial \theta_i}{\partial u_j}-\frac{\partial \theta_i}{\partial u_k}<0$
for $(r_i, r_j, r_k)\in \Omega^E_{ijk}(\eta)\cap U$.
By the explicit form of $V_i$, i.e.
\begin{equation*}
\begin{aligned}
V_i=\{(r_i, r_j, r_k)\in \mathbb{R}^3_{>0}| \kappa_i\geq \frac{-B_i+\sqrt{\Delta_i}}{2A_i}\}
   =\{(r_i, r_j, r_k)\in \mathbb{R}^3_{>0}| r_i\leq \frac{2A_i}{-B_i+\sqrt{\Delta_i}}\},
\end{aligned}
\end{equation*}
we have $\theta_i\rightarrow \pi$ as $(r_i, r_j, r_k)\rightarrow (\overline{r}_i, \overline{r}_j, \overline{r}_k)$.
Otherwise, $\theta_i\rightarrow 0$ as $(r_i, r_j, r_k)\rightarrow (\overline{r}_i, \overline{r}_j, \overline{r}_k)$.
As a result,  by $\frac{\partial \theta_i}{\partial u_i}<0$, we have $\theta_i< 0$ for
$(\overline{r}_i+\epsilon, \overline{r}_j, \overline{r}_k)\in \Omega^E_{ijk}(\eta)\cap U$, $\epsilon>0$ small enough. It is impossible.
By $\theta_i+\theta_j+\theta_k=\pi$, we have $\theta_j\rightarrow 0$, $\theta_k\rightarrow 0$ as
$(r_i, r_j, r_k)\rightarrow (\overline{r}_i, \overline{r}_j, \overline{r}_k)$.
The same arguments apply to the other components of $\mathbb{R}^3_{>0}\setminus\Omega_{ijk}^E(\eta)$.

Therefore, the extension (\ref{extension of theta_i}) defines a continuous extension of  the inner angle functions $\theta_i, \theta_j, \theta_k$
 on $\mathbb{R}^3_{>0}$.
 \qed

By Lemma \ref{Euclidean extension}, we can extend the combinatorial curvature function $K$ defined for nondegenerate
radius vectors to be defined for all $r\in \mathbb{R}^N_{>0}$ by setting
\begin{equation}\label{extension of combinatorial curvature}
\begin{aligned}
\widetilde{K}_i=2\pi-\sum_{\{ijk\}\in F}\widetilde{\theta}_i,
\end{aligned}
\end{equation}
where $\widetilde{\theta}_i$ is the extension of $\theta_i$ defined  by (\ref{extension of theta_i}).
The extended combinatorial curvature $\widetilde{K}$ still satisfies the discrete Gauss-Bonnet formula
$\sum_{i=1}^N\widetilde{K}_i=2\pi\chi(M)$.

Recall the following definition of closed continuous $1$-form and
extension of locally convex function of Luo \cite{L4}, which is a generalization
of Bobenko-Pinkall-Spingborn's extension introduced in \cite{BPS}.
\begin{definition}[\cite{L4}, Definition 2.3]
A differential 1-form $w=\sum_{i=1}^n a_i(x)dx^i$ in an open set $U\subset \mathbb{R}^n$ is said to be continuous if each $a_i(x)$ is continuous on $U$. A continuous differential 1-form $w$ is called closed if $\int_{\partial \tau}w=0$ for each
triangle $\tau\subset U$.
\end{definition}

\begin{theorem}[\cite{L4}, Corollary 2.6]\label{Luo's convex extention}
Suppose $X\subset \mathbb{R}^n$ is an open convex set and $A\subset X$ is an open subset of $X$ bounded by a real analytic codimension-1 submanifold in $X$. If $w=\sum_{i=1}^na_i(x)dx_i$ is a continuous closed 1-form on $A$ so that $F(x)=\int_a^x w$ is locally convex on $A$ and each $a_i$ can be extended continuous to $X$ by constant functions to a function $\widetilde{a}_i$ on $X$, then  $\widetilde{F}(x)=\int_a^x\sum_{i=1}^n\widetilde{a}_i(x)dx_i$ is a $C^1$-smooth
convex function on $X$ extending $F$.
\end{theorem}

By Lemma \ref{Euclidean extension} and Theorem \ref{Luo's convex extention}, the locally concave function
$\mathcal{E}_{ijk}$ defined by (\ref{Euclidean Ricci energy function for triangle})
for nondegenerate $(u_i,u_j,u_k)$ on $(\{ijk\}, \varepsilon,\eta)$
can be extended to be a $C^1$ smooth concave function
\begin{equation}\label{extension of Euclidean Ricci energy function for triangle}
\begin{aligned}
\widetilde{\mathcal{E}}_{ijk}(u_i,u_j,u_k)=\int_{(\overline{u}_i, \overline{u}_j, \overline{u}_k)}^{(u_i, u_j, u_k)}\widetilde{\theta}_idu_i+\widetilde{\theta}_jdu_j+\widetilde{\theta}_kdu_k
\end{aligned}
\end{equation}
defined for all $(u_i,u_j,u_k)\in \mathbb{R}^3$ with $\nabla_{u_i}\widetilde{\mathcal{E}}_{ijk}=\widetilde{\theta}_i$.
As a result, the locally convex function $\mathcal{E}$ defined by (\ref{Euclidean Ricci energy function}) for nondegenerate Euclidean
discrete conformal factors can be extended to
be a $C^1$ smooth convex function
\begin{equation}\label{extended Euclidean energy function}
\begin{aligned}
\widetilde{\mathcal{E}}(u_1, \cdots, u_N)=2\pi\sum_{i\in V}u_i-\sum_{\{ijk\}\in F}\widetilde{\mathcal{E}}_{ijk}(u_i,u_j,u_k)
\end{aligned}
\end{equation}
defined on $\mathbb{R}^N$ with $\nabla_{u_i} \widetilde{\mathcal{E}}=\widetilde{K}_i=2\pi-\sum\widetilde{\theta}_i$.

Using the extended Ricci energy function $\widetilde{\mathcal{E}}$, we can prove the following rigidity for Glickenstein's
Euclidean discrete conformal structures on polyhedral surfaces,
which is a generalization of Theorem \ref{main rigidity introduction} (a).

\begin{theorem}\label{Thm rigidity of EDCS context}
Suppose $(M, \mathcal{T}, \varepsilon, \eta)$ is a weighted triangulated surface with
the weights $\varepsilon: V\rightarrow \{0, 1\}$ and $\eta: E\rightarrow \mathbb{R}$
satisfying the structure conditions (\ref{structure condition 1}) and (\ref{structure condition 2}).
If there exists a nondegenrate radius vector $r_A\in \Omega^E$ and
a radius vector $r_B\in \mathbb{R}^N_{>0}$ such that
$K(r_A)=\widetilde{K}(r_B)$. Then $r_A=cr_B$ for some positive constant $c\in \mathbb{R}$.
\end{theorem}

\proof
Set
\begin{equation*}
\begin{aligned}
\mathcal{F}(t)=\widetilde{\mathcal{E}}((1-t)u_A+tu_B)=2\pi\sum_{i=1}^N[(1-t)u_{A,i}+tu_{B,i}]+\sum_{\{ijk\}\in F}\mathcal{F}_{ijk}(t),
\end{aligned}
\end{equation*}
where
$\mathcal{F}_{ijk}(t)=-\widetilde{\mathcal{E}}_{ijk}((1-t)u_A+tu_B)$.
Then $\mathcal{F}(t)$ is a $C^1$ smooth convex function for $t\in [0, 1]$ with $\mathcal{F}'(0)=\mathcal{F}'(1)$.
This implies $\mathcal{F}'(t)=\mathcal{F}'(0)$ for any $t\in [0, 1]$.
Note that the admissible space $\Omega^E$ of nondegenerate Euclidean discrete conformal factors is an
open subset of $\mathbb{R}^N$, there exists $\epsilon>0$ such that $(1-t)u_A+tu_B$ is nondegenerate for $t\in [0, \epsilon]$.
Note that $\mathcal{F}(t)$ is smooth for $t\in [0, \epsilon]$,  by $\mathcal{F}'(t)=\mathcal{F}'(0)$ for any $t\in [0, 1]$,  we have
\begin{equation*}
\begin{aligned}
\mathcal{F}''(t)=(u_B-u_A) \Lambda^E (u_B-u_A)^T=0, \ \forall t\in [0, \epsilon].
\end{aligned}
\end{equation*}
By  Corollary \ref{Euclidean curvature Jacobian positivity}, this implies $u_B-u_A=\lambda(1, \cdots, 1)$
for some constant $\lambda\in \mathbb{R}$.
As a result, $r_B=cr_A$ with $c=e^\lambda>0$.
\qed

\section{Hyperbolic discrete conformal structures}\label{section 3}
\subsection{Admissible space of hyperbolic discrete conformal factors for a triangle}
In this subsection,
we fix a weighted triangle $(\sigma=\{ijk\}, \varepsilon,\eta)$ with two weights
$\varepsilon: V_\sigma\rightarrow \{0, 1\}$ and $E_\sigma\rightarrow \mathbb{R}$
satisfying the structure conditions (\ref{structure condition 1}) and (\ref{structure condition 2}).
In the hyperbolic background geometry,
the lengths $l_{ij}, l_{ik}, l_{jk}$ of the edges in $E_\sigma$ are defined by the discrete conformal factor $f: V_\sigma \rightarrow \mathbb{R}$ via the formula (\ref{defn of hyperbolic length}).
The discrete conformal factor $f$ is \textit{nondegenerate} if $l_{ij}, l_{ik}, l_{jk}$ satisfy the triangle inequality, otherwise it is \textit{degenerate}.
We use $\Omega_{ijk}^H(\eta)$ to denote the space of nondegenerate hyperbolic discrete conformal factors for the weighted triangle $(\{ijk\}, \varepsilon,\eta)$.
In this subsection, we will give an analytical characterization of $\Omega_{ijk}^H(\eta)$.
The method is a modification of the Euclidean case.
As many results in this subsection are paralleling to the results in Subsection \ref{subsection Euclidean admissible space},
some proofs for the results in this subsection will be omitted if there is no difference.

To simplify the notations, set
\begin{equation}\label{simplification C_i S_i}
\begin{aligned}
S_i=e^{f_i}, C_i=\sqrt{1+\varepsilon_ie^{2f_i}}, \kappa_i=\frac{C_i}{S_i}.
\end{aligned}
\end{equation}
Then
\begin{equation}\label{relation of C_i S_i}
\begin{aligned}
C_i^2-\varepsilon_iS_i^2=1
\end{aligned}
\end{equation}
and the hyperbolic edge length $l_{ij}$ is determined by
\begin{equation}\label{hyperbolic edge length in C_i S_i}
\begin{aligned}
\cosh l_{ij}=C_iC_j+\eta_{ij}S_iS_j.
\end{aligned}
\end{equation}

By the structure condition  (\ref{structure condition 1})  and  the inequality
$(1+a^2)(1+b^2)\geq (1+ab)^2$,
we have
$$\sqrt{(1+\varepsilon_ie^{2f_i})(1+\varepsilon_je^{2f_j})}+\eta_{ij}e^{f_i+f_j}
\geq 1+(\varepsilon_i\varepsilon_j+\eta_{ij})e^{f_i+f_j}> 1.$$
This implies the hyperbolic edge length defined by (\ref{defn of hyperbolic length}) is well defined.
Parallelling to Lemma \ref{lemma triangle equality euqi to Q>0},
we have the following result on the triangle inequalities in the hyperbolic background geometry.

\begin{lemma}\label{lemma hyper triangle ineq in Q}
For the weighted triangle $(\{ijk\}, \varepsilon,\eta)$,  the edge lengths
$l_{ij}, l_{ik}, l_{jk}$ defined by (\ref{defn of hyperbolic length})
satisfy the triangle inequalities if and only if $Q^H>0$,
where
\begin{equation}\label{hyperbolic Q}
\begin{aligned}
Q^H=&(\varepsilon_j\varepsilon_k-\eta_{jk}^2)\kappa_i^2+(\varepsilon_i\varepsilon_k-\eta_{ik}^2)\kappa_j^2+(\varepsilon_i\varepsilon_j-\eta_{ij}^2)\kappa_k^2
 +2\gamma_i\kappa_j\kappa_k+2\gamma_j\kappa_i\kappa_k+2\gamma_k\kappa_i\kappa_j\\
 &
 +\varepsilon_i\eta_{jk}^2
  +\varepsilon_j\eta_{ik}^2+\varepsilon_k\eta_{ij}^2+2\eta_{ij}\eta_{ik}\eta_{jk}-\varepsilon_i\varepsilon_j\varepsilon_k.
\end{aligned}
\end{equation}
\end{lemma}
\proof
Note that the positive edge lengths $l_{ij}, l_{ik}, l_{jk}$ defined by (\ref{defn of hyperbolic length})
satisfy the triangle inequalities if and only if
%
\begin{equation}\label{equivalent hyper triangle inequ proof}
\begin{aligned}
0<&4\sinh\frac{l_{ij}+l_{ik}+l_{jk}}{2}\sinh\frac{l_{ij}+l_{ik}-l_{jk}}{2}\sinh\frac{l_{ij}-l_{ik}+l_{jk}}{2}\sinh\frac{-l_{ij}+l_{ik}+l_{jk}}{2}\\
=&(\cosh (l_{ik}+l_{jk})-\cosh l_{ij})(\cosh l_{ij}-\cosh (l_{jk}-l_{ik}))\\
=&1+2\cosh l_{ij}\cosh l_{ik}\cosh l_{jk}-\cosh^2l_{ij}-\cosh^2l_{ik}-\cosh^2l_{jk}.
\end{aligned}
\end{equation}
Submitting (\ref{hyperbolic edge length in C_i S_i}) into (\ref{equivalent hyper triangle inequ proof}) and by direct calculations, we have
\begin{equation}\label{expansion of hyperbolic triangle inequalities}
\begin{aligned}
0<&4\sinh\frac{l_{ij}+l_{ik}+l_{jk}}{2}\sinh\frac{l_{ij}+l_{ik}-l_{jk}}{2}\sinh\frac{l_{ij}-l_{ik}+l_{jk}}{2}\sinh\frac{-l_{ij}+l_{ik}+l_{jk}}{2}\\
=&(1+2C_i^2C_j^2C_k^2-C_i^2C_j^2-C_i^2C_k^2-C_j^2C_k^2)\\
 &+2\gamma_iC_jC_kS_i^2S_jS_k+2\gamma_jC_iC_kS_iS_j^2S_k+2\gamma_kC_iC_jS_iS_jS_k^2\\
 &+2\eta_{ij}\eta_{ik}\eta_{jk}S_i^2S_j^2S_k^2-\eta_{ij}^2S_i^2S_j^2-\eta_{ik}^2S_i^2S_k^2-\eta_{jk}^2S_j^2S_k^2\\
=&(\varepsilon_i\varepsilon_j-\eta_{ij}^2)S_i^2S_j^2+(\varepsilon_i\varepsilon_k-\eta_{ik}^2)S_i^2S_k^2
  +(\varepsilon_j\varepsilon_k-\eta_{jk}^2)S_j^2S_k^2
 +2\gamma_iC_jC_kS_i^2S_jS_k\\
 &+2\gamma_jC_iC_kS_iS_j^2S_k+2\gamma_kC_iC_jS_iS_jS_k^2
 +2(\varepsilon_i\varepsilon_j\varepsilon_k+\eta_{ij}\eta_{ik}\eta_{jk})S_i^2S_j^2S_k^2\\
=&S_i^2S_j^2S_k^2
  [(\varepsilon_j\varepsilon_k-\eta_{jk}^2)\kappa_i^2+(\varepsilon_i\varepsilon_k-\eta_{ik}^2)\kappa_j^2+(\varepsilon_i\varepsilon_j-\eta_{ij}^2)\kappa_k^2\\
 &+2\gamma_i\kappa_j\kappa_k+2\gamma_j\kappa_i\kappa_k+2\gamma_k\kappa_i\kappa_j
 +\varepsilon_i\eta_{jk}^2
  +\varepsilon_j\eta_{ik}^2+\varepsilon_k\eta_{ij}^2+2\eta_{ij}\eta_{ik}\eta_{jk}-\varepsilon_i\varepsilon_j\varepsilon_k]\\
=&S_i^2S_j^2S_k^2Q^H,
\end{aligned}
\end{equation}
where (\ref{relation of C_i S_i}) is used in the second and third equality.
This completes the proof.
\qed

Comparing Lemma \ref{lemma triangle equality euqi to Q>0} with Lemma \ref{lemma hyper triangle ineq in Q},
we find that $Q^H$  in Lemma \ref{lemma hyper triangle ineq in Q}
has one term $G=\varepsilon_i\eta_{jk}^2
  +\varepsilon_j\eta_{ik}^2+\varepsilon_k\eta_{ij}^2+2\eta_{ij}\eta_{ik}\eta_{jk}-\varepsilon_i\varepsilon_j\varepsilon_k$ more than
$Q^E$  in Lemma \ref{lemma triangle equality euqi to Q>0}.
Furthermore, this term is symmetric in $i,j,k$.
Set $h_i, h_j, h_k$ as that in (\ref{hi hj hk}), then
\begin{equation}\label{relation of hyperbolic Q and hi hj hk}
\begin{aligned}
Q^H=\kappa_ih_i+\kappa_jh_j+\kappa_kh_k+G=Q^E+G.
\end{aligned}
\end{equation}

Paralleling to Theorem \ref{theorem admissible space Euclidean}, we have the following
analytic characterization of the admissible space $\Omega^H_{ijk}(\eta)$
of hyperbolic discrete conformal factors for the weighted triangle $(\{ijk\}, \varepsilon,\eta)$.

\begin{theorem}\label{hyperbolic admissible space}
For the weighted triangle $(\{ijk\}, \varepsilon,\eta)$,
the admissible space $\Omega_{ijk}^H(\eta)\subseteq \mathbb{R}^3$
is a nonempty simply connected open set whose boundary components are analytical.
Furthermore,
\begin{equation*}
\begin{aligned}
\Omega^H_{ijk}(\eta)=\mathbb{R}^3\setminus \sqcup_{\alpha\in \Lambda} V_\alpha,
\end{aligned}
\end{equation*}
where $\Lambda=\{q\in \{i,j,k\}|A_q=\eta_{st}^2-\varepsilon_s\varepsilon_t>0, \{q, s, t\}=\{i, j, k\}\}$,
$\sqcup_{\alpha\in \Lambda}V_\alpha$ is a disjoint union of $V_\alpha$ with
\begin{equation}\label{discription of V_i hyperbolic}
\begin{aligned}
V_i=&\left\{(f_i, f_j, f_k)\in \mathbb{R}^3| \kappa_i\geq \frac{-B_i+\sqrt{\Delta_i}}{2A_i}\right\}\\
   =&\left\{(f_i, f_j, f_k)\in \mathbb{R}^3| f_i\leq -\frac{1}{2} \ln\left[\left(\frac{-B_i+\sqrt{\Delta_i}}{2A_i}\right)^2-\varepsilon_i\right]\right\}
\end{aligned}
\end{equation}
being a closed region in $\mathbb{R}^3$ bounded by an analytical function defined on $\mathbb{R}^2$
and $V_j, V_k$ defined similarly.
Here $\Delta_i=B_i^2-4A_iC_i$ and $A_i, B_i, C_i$ are defined by (\ref{hyperbolic Ai Bi Ci}).
\end{theorem}

To prove Theorem \ref{hyperbolic admissible space}, we first prove the following result for the hyperbolic discrete conformal structures
paralleling to Lemma \ref{A_i nonpositive} for the Euclidean discrete conformal structures.

\begin{lemma}\label{hyperbolic A_i nonpositive implies simply connected}
For the weighted triangle $(\{ijk\}, \varepsilon,\eta)$,
if $\varepsilon_j\varepsilon_k-\eta_{jk}^2\geq 0$, $\varepsilon_i\varepsilon_k-\eta_{ik}^2\geq 0$ and
$\varepsilon_i\varepsilon_j-\eta_{ij}^2\geq 0$,
then the admissible space $\Omega^H_{ijk}(\eta)$ of hyperbolic discrete conformal factors $(f_i,f_j,f_k)$ is $\mathbb{R}^3$
and hence simply connected.
\end{lemma}
\proof
By Lemma \ref{lemma hyper triangle ineq in Q}, we just need to prove $Q^H>0$ for any $(f_i,f_j,f_k)\in \mathbb{R}^3$.

If one of $\varepsilon_i, \varepsilon_j, \varepsilon_k$ is zero,
say $\varepsilon_i=0$, we have $\eta_{ij}>0$, $\eta_{ik}>0$, $\eta_{jk}+\varepsilon_j\varepsilon_k$ by the structure condition (\ref{structure condition 1}).
This implies $G=2\eta_{ij}\eta_{ik}\eta_{jk}+\varepsilon_j\eta_{ik}^2+\varepsilon_k\eta_{ij}^2
=(\varepsilon_j\eta_{ik}-\varepsilon_k\eta_{ij})^2+2\eta_{ij}\eta_{ik}(\eta_{jk}+\varepsilon_j\varepsilon_k)>0$.
By the structure condition (\ref{structure condition 2}) and the condition $\varepsilon_j\varepsilon_k-\eta_{jk}^2\geq 0$, $\varepsilon_i\varepsilon_k-\eta_{ik}^2\geq 0$, $\varepsilon_i\varepsilon_j-\eta_{ij}^2\geq 0$,
we have $h_i\geq 0, h_j\geq 0, h_k\geq 0$.
By (\ref{relation of hyperbolic Q and hi hj hk}), this implies $Q^H\geq G>0$.

If $\varepsilon_i=\varepsilon_j=\varepsilon_k=1$, by (\ref{simplification C_i S_i}), we have $\kappa_i>1, \kappa_j>1, \kappa_j>1$.
Combining the structure condition (\ref{structure condition 2}) and the assumption $\varepsilon_j\varepsilon_k-\eta_{jk}^2\geq 0$, $\varepsilon_i\varepsilon_k-\eta_{ik}^2\geq 0$, $\varepsilon_i\varepsilon_j-\eta_{ij}^2\geq 0$, this implies
\begin{equation*}
\begin{aligned}
Q^H\geq&1-\eta_{jk}^2+1-\eta_{ik}^2+1-\eta_{ij}^2
 +2\gamma_i+2\gamma_j+2\gamma_k
 +\eta_{jk}^2+\eta_{ik}^2+\eta_{ij}^2+2\eta_{ij}\eta_{ik}\eta_{jk}-1\\
 =&2(\eta_{ij}+1)(\eta_{ik}+1)(\eta_{jk}+1)\\
 >&0,
\end{aligned}
\end{equation*}
where the structure condition (\ref{structure condition 1})
is used in the last inequality.
Therefore, the admissible space $\Omega^H_{ijk}(\eta)$ is $\mathbb{R}^3$.
\qed

By Lemma \ref{hyperbolic A_i nonpositive implies simply connected}, we just need to study the admissible space $\Omega^H_{ijk}(H)$ for
the case that one of $\varepsilon_j\varepsilon_k-\eta_{jk}^2$, $\varepsilon_i\varepsilon_k-\eta_{ik}^2$, $\varepsilon_i\varepsilon_j-\eta_{ij}^2$
is negative.

Parallelling to Lemma \ref{sign of hi hj hk}, we have the following result on the signs of $h_i, h_j, h_k$ for the degenerate hyperbolic discrete conformal factors on $(\{ijk\}, \varepsilon,\eta)$.

\begin{lemma}\label{hyperbolic signs of hi hj hk}
For the weighted triangle $(\{ijk\}, \varepsilon,\eta)$,
if $(f_i, f_j, f_k)$ is a degenerate hyperbolic discrete conformal factor,
then one of $h_i, h_j, h_k$ is negative and the other two are positive.
\end{lemma}
\proof
By Lemma \ref{lemma hyper triangle ineq in Q}, if $(f_i, f_j, f_k)$ is a degenerate hyperbolic discrete conformal factor
for $(\{ijk\}, \varepsilon,\eta)$, then
$Q^H=\kappa_ih_i+\kappa_jh_j+\kappa_kh_k+G\leq 0$.
By Lemma \ref{hyperbolic A_i nonpositive implies simply connected}, at least one of
$\varepsilon_j\varepsilon_k-\eta_{jk}^2$, $\varepsilon_i\varepsilon_k-\eta_{ik}^2$, $\varepsilon_i\varepsilon_j-\eta_{ij}^2$
is negative.
Without loss of generality, assume
$\varepsilon_j\varepsilon_k-\eta_{jk}^2<0$.
By Corollary \ref{G positive}, we have $G>0$.
By $Q^H=\kappa_ih_i+\kappa_jh_j+\kappa_kh_k+G\leq0$, this implies
$Q^E=\kappa_ih_i+\kappa_jh_j+\kappa_kh_k\leq -G< 0$.
Therefore,  at least one of $h_i, h_j, h_k$ is negative.
Following the proof for Lemma \ref{sign of hi hj hk}, we have one of $h_i, h_j, h_k$ is negative and the other two are positive.
As the proof is parallelling to that for Lemma \ref{sign of hi hj hk}, we omit the details here.
\qed

Different from the Euclidean case, we need to use Lemma \ref{hyperbolic A_i nonpositive implies simply connected}
to prove Lemma \ref{hyperbolic signs of hi hj hk} in the hyperbolic case.

Suppose $(f_i,f_j,f_k)\in \mathbb{R}^3$ is a degenerate hyperbolic discrete conformal factor on $(\{ijk\}, \varepsilon,\eta)$.
By Lemma \ref{hyperbolic signs of hi hj hk}, one of $h_i, h_j, h_k$ is negative.
Without loss of generality, assume $h_i<0$ at $(f_i,f_j,f_k)$.
By the structure condition (\ref{structure condition 2}),
this implies $(\eta_{jk}^2-\varepsilon_j\varepsilon_k)\kappa_i>\gamma_j\kappa_k+\gamma_k\kappa_j\geq 0$.
As $(f_i,f_j,f_k)$ is a degenerate hyperbolic discrete conformal factor,
we have $Q^H\leq 0$ by Lemma  \ref{lemma hyper triangle ineq in Q}.
This is equivalent to
\begin{equation}\label{hyperbolic quadratic inequality}
\begin{aligned}
A_i\kappa_i^2+B_i\kappa_i+C_i\geq 0,
\end{aligned}
\end{equation}
where
\begin{equation}\label{hyperbolic Ai Bi Ci}
\begin{aligned}
A_i=&\eta_{jk}^2-\varepsilon_j\varepsilon_k>0,\\
B_i=&-2\gamma_j\kappa_k-2\gamma_k\kappa_j\leq 0,\\
C_i=&(\eta_{ik}^2-\varepsilon_i\varepsilon_k)\kappa_j^2+(\eta_{ij}^2-\varepsilon_i\varepsilon_j)\kappa_k^2-2\gamma_i\kappa_j\kappa_k-G.
\end{aligned}
\end{equation}
Parallelling to Lemma \ref{discriminant of kappa i}, we have the following result for
the discriminant of (\ref{hyperbolic quadratic inequality}) in the hyperbolic case.
\begin{lemma}\label{hyperbolic nonnegative discriminant}
For the weighted triangle $(\{ijk\}, \varepsilon,\eta)$,
if $A_i=\eta_{jk}^2-\varepsilon_j\varepsilon_k>0$, then the discriminant $\Delta_i=B_i^2-4A_iC_i$ for (\ref{hyperbolic quadratic inequality})
is positive,
where $A_i, B_i, C_i$ are defined by (\ref{hyperbolic Ai Bi Ci}).
\end{lemma}
\proof
By the assumption $A_i=\eta_{jk}^2-\varepsilon_j\varepsilon_k>0$ and Corollary \ref{G positive}, we have $G> 0$.
Then the proof is reduced to the case in
Lemma \ref{discriminant of kappa i}, which has been completed.
\qed

\begin{remark}\label{discriminant of kappa j k hyperbolic}
One can also take $Q^H$ as a quadratic function of $\kappa_j$ or $\kappa_k$ and define $\Delta_j$, $\Delta_k$
similarly.
By symmetry, we have
$\Delta_j> 0$ if $A_j=\eta_{ik}^2-\varepsilon_i\varepsilon_k>0$ and
$\Delta_k>0$ if $A_k=\eta_{ij}^2-\varepsilon_i\varepsilon_j>0$.
\end{remark}

Note that we have Lemma \ref{lemma hyper triangle ineq in Q}, Lemma \ref{hyperbolic A_i nonpositive implies simply connected},
Lemma \ref{hyperbolic signs of hi hj hk}, Lemma \ref{hyperbolic nonnegative discriminant} in the hyperbolic case,
which are paralleling to
Lemma \ref{lemma triangle equality euqi to Q>0}, Lemma \ref{A_i nonpositive}, Lemma \ref{sign of hi hj hk},
Lemma \ref{discriminant of kappa i} in the Euclidean case respectively.
Then the proof of Theorem \ref{hyperbolic admissible space} is parallelling to that of Theorem \ref{theorem admissible space Euclidean}.
We omit the details here.

\begin{remark}\label{hi positive in Vi hyperbolic}
Suppose $(f_i, f_j, f_k)\in V_i$ is a degenerate hyperbolic discrete conformal factor for the weighted triangle $(\{ijk\}, \varepsilon,\eta)$.
Then by $-2h_i=2A_i\kappa_i+B_i$ and (\ref{discription of V_i hyperbolic}), we have $h_i< 0$.
By Lemma \ref{hyperbolic signs of hi hj hk}, this implies $h_i<0, h_j>0, h_k>0$.
\end{remark}

Define
\begin{equation*}
\begin{aligned}
\Omega_{ijk}^H=\left\{(f_i, f_j, f_k, \eta_{ij}, \eta_{ik}, \eta_{jk})\in \mathbb{R}^6|
\eta \ \text{satisfies}\ (\ref{structure condition 1}), (\ref{structure condition 2}) \ \text{and } (f_i, f_j, f_k)\in \Omega_{ijk}^H(\eta)\right\}.
\end{aligned}
\end{equation*}
As a corollary of Theorem \ref{hyperbolic admissible space}, we have the following result for
the parameterized hyperbolic admissible space $\Omega_{ijk}^H$.
\begin{corollary}\label{coro connectivity of para admi space hyperbolic}
For the triangle $\sigma=\{ijk\}$ with a weight $\varepsilon: V_\sigma\rightarrow \{0, 1\}$,  the parameterized hyperbolic admissible space $\Omega_{ijk}^H$ is connected.
\end{corollary}

The proof for Corollary \ref{coro connectivity of para admi space hyperbolic}
is the same as that for Corollary \ref{coro connectivity of para admi space}, so
we omit the details of the proof here.
Parallelling to the Euclidean case, the parameterized admissible space $\Omega_{ijk}^H$
contains some points with good properties.

\begin{lemma}\label{hi sign at special point hyperbolic}
The point $p=(f_i, f_j, f_k, \eta_{ij}, \eta_{ik}, \eta_{jk})=(0, 0, 0, 1, 1, 1)$ is a point in $\Omega_{ijk}^H$.
Furthermore,
$h_i(p)>0, h_j(p)>0, h_k(p)>0$.
\end{lemma}
\proof
It is straight forward to check that $(\eta_{ij}, \eta_{ik}, \eta_{jk})=(1,1,1)$
satisfies the structure conditions (\ref{structure condition 1}) and (\ref{structure condition 2}).
By $\varepsilon_i,\varepsilon_j,\varepsilon_k\in \{0,1\}$, we have
\begin{equation*}
\begin{aligned}
h_i(p)=&(\varepsilon_j\varepsilon_k-1)\sqrt{1+\varepsilon_i}+(1+\varepsilon_j)\sqrt{1+\varepsilon_k}+(1+\varepsilon_k)\sqrt{1+\varepsilon_j}\\
\geq &\sqrt{1+\varepsilon_j}+\sqrt{1+\varepsilon_k}-\sqrt{1+\varepsilon_i}\\
\geq &2-\sqrt{2}\\
>&0.
\end{aligned}
\end{equation*}
Similarly, we have $h_j(p)>0, h_k(p)>0$.
On the other hand, by $\varepsilon_i,\varepsilon_j,\varepsilon_k\in \{0,1\}$, we have
$G(p)=2+\varepsilon_i+\varepsilon_j+\varepsilon_k-\varepsilon_i\varepsilon_j\varepsilon_k>0$.
As a result,
we have $Q^H=\kappa_ih_i+\kappa_jh_j+\kappa_kh_k+G>0$
at $p$. Therefore, $p\in \Omega_{ijk}^H$.
\qed

\subsection{Negative definiteness of the Jacobian matrix in the hyperbolic background geometry}
Let $(\{ijk\}, \varepsilon,\eta)$ be a nondegenerate hyperbolic weighted triangle
with edge lengths given by (\ref{defn of hyperbolic length}). Suppose
$\theta_i, \theta_j, \theta_k$ are the inner angles at the vertices $i, j, k$ in the triangle respectively.
Set
\begin{equation}\label{ui in fi hyperbolic}
\begin{aligned}
u_i=
\left\{
  \begin{array}{ll}
    f_i, & \hbox{$\varepsilon_i=0$;} \\
    \frac{1}{2}\ln\left(\frac{\sqrt{1+e^{2f_i}}-1}{\sqrt{1+e^{2f_i}}+1}\right), & \hbox{$\varepsilon_i=1$.}
  \end{array}
\right.
\end{aligned}
\end{equation}
Then
\begin{equation}\label{derivative of fi in ui hyperbolic}
\begin{aligned}
\frac{\partial f_i}{\partial u_i}=\sqrt{1+\varepsilon_ie^{2f_i}}=C_i.
\end{aligned}
\end{equation}

\begin{lemma}[\cite{GT}]\label{hyperbolic symmetry}
Let $(\{ijk\}, \varepsilon,\eta)$ be a weighted triangle and $(f_i, f_j, f_k)\in \mathbb{R}^3$ is a nondegenerate hyperbolic
discrete conformal factor on $(\{ijk\}, \varepsilon,\eta)$.
Then
\begin{equation}\label{derivative of theta i of uj hyperbolic}
\begin{aligned}
\frac{\partial \theta_i}{\partial u_j}
=\frac{\partial \theta_j}{\partial u_i}
=\frac{S_i^2S_j^2S_k}{A\sinh^2 l_{ij}}[(\varepsilon_i\varepsilon_j-\eta_{ij}^2)\kappa_k+\gamma_i\kappa_j+\gamma_j\kappa_i]
=\frac{S_i^2S_j^2S_k}{A\sinh^2 l_{ij}}h_k,
\end{aligned}
\end{equation}
where $A=\sinh l_{ij}\sinh l_{ik}\sin \theta_i$ and $u_i$ is defined by (\ref{ui in fi hyperbolic}).
\end{lemma}
\proof
By the chain rules,
\begin{equation}\label{chain rules for hyperbolic triangle}
\begin{aligned}
\frac{\partial \theta_i}{\partial u_j}
=\frac{\partial \theta_i}{\partial l_{jk}}\frac{\partial l_{jk}}{\partial u_j}
 +\frac{\partial \theta_i}{\partial l_{ik}}\frac{\partial l_{ik}}{\partial u_j}
 +\frac{\partial \theta_i}{\partial l_{ij}}\frac{\partial l_{ij}}{\partial u_j}.
\end{aligned}
\end{equation}
The derivative cosine law (\cite{CL}, Lemma A1) for hyperbolic triangles gives
\begin{equation}\label{hyperbolic derivative cosine law}
\begin{aligned}
\frac{\partial \theta_i}{\partial l_{jk}}=\frac{\sinh l_{jk}}{A},\
\frac{\partial \theta_i}{\partial l_{ij}}=\frac{-\sinh l_{jk}\cos\theta_j}{A},
\end{aligned}
\end{equation}
where $A=\sinh l_{ij}\sinh l_{ik}\sin\theta_i$.
By (\ref{defn of hyperbolic length}) and (\ref{derivative of fi in ui hyperbolic}), we have
\begin{equation}\label{derivative of hyperbolic length}
\begin{aligned}
\frac{\partial l_{jk}}{\partial u_j}=\frac{1}{\sinh l_{jk}} (\varepsilon_jS_j^2C_k+\eta_{jk}S_jS_kC_j),\
\frac{\partial l_{ik}}{\partial u_j}=0,\
\frac{\partial l_{ij}}{\partial u_j}=\frac{1}{\sinh l_{ij}}(\varepsilon_jS_j^2C_i+\eta_{ij}S_iS_jC_j).
\end{aligned}
\end{equation}
Submitting (\ref{hyperbolic derivative cosine law}) and (\ref{derivative of hyperbolic length})
into (\ref{chain rules for hyperbolic triangle}), by direct calculations, we have
\begin{equation}\label{hyperbolic derivative of thetai in uj proof}
\begin{aligned}
\frac{\partial \theta_i}{\partial u_j}
=&\frac{1}{A}  (\varepsilon_jS_j^2C_k+\eta_{jk}S_jS_kC_j)
 +\frac{-\sinh l_{jk}\cos\theta_j}{A}\frac{1}{\sinh l_{ij}}(\varepsilon_jS_j^2C_i+\eta_{ij}S_iS_jC_j)\\
=&\frac{1}{A\sinh^2 l_{ij}}[(\cosh^2l_{ij}-1)(\varepsilon_jS_j^2C_k+\eta_{jk}S_jS_kC_j)\\
  &+(\cosh l_{ik}-\cosh l_{ij}\cosh l_{jk})(\varepsilon_jS_j^2C_i+\eta_{ij}S_iS_jC_j)]\\
=&\frac{S_i^2S_j^2S_k}{A\sinh^2 l_{ij}}[(\varepsilon_i\varepsilon_j-\eta_{ij}^2)\kappa_k+\gamma_i\kappa_j+\gamma_j\kappa_i]\\
=&\frac{S_i^2S_j^2S_k}{A\sinh^2 l_{ij}}h_k,
\end{aligned}
\end{equation}
where the hyperbolic cosine law is used in the second equality and
the definition (\ref{defn of hyperbolic length}) for hyperbolic length is used in the third equality.
Note that (\ref{hyperbolic derivative of thetai in uj proof}) is symmetric in the indices $i$ and $j$, we have
$\frac{\partial \theta_i}{\partial u_j}=\frac{\partial \theta_j}{\partial u_i}$.
\qed

\begin{remark}\label{derivative tends infty hyperbolic}
The result
in Lemma \ref{hyperbolic symmetry} was proved by Glickenstein-Thomas \cite{GT} and Zhang-Guo-Zeng-Luo-Yau-Gu \cite{ZGZLYG}.
Here we give a proof by direct calculations for completeness.
By (\ref{derivative of theta i of uj hyperbolic}) and Remark \ref{hi positive in Vi hyperbolic},
if $(f_i, f_j, f_k)\in \Omega_{ijk}^H(\eta)$ tends to a point $(\overline{f}_i, \overline{f}_j, \overline{f}_k)\in \partial V_i$ with $V_i\neq \emptyset$,
 then $\frac{\partial \theta_i}{\partial u_j}\rightarrow +\infty$, $\frac{\partial \theta_i}{\partial u_k}\rightarrow +\infty$.
Recall the following formula obtained by Glickenstein-Thomas (\cite{GT}, Proposition 9)
\begin{equation}\label{GT's formula}
\frac{\partial A_{ijk}}{\partial u_i}=\frac{\partial \theta_j}{\partial u_i}(\cosh l_{ij}-1)+\frac{\partial \theta_k}{\partial u_i}(\cosh l_{ik}-1)
\end{equation}
for the area $A_{ijk}$ of the hyperbolic triangle $\{ijk\}$,
we have
\begin{equation}\label{relation of partial derivative hyperbolic}
-\frac{\partial \theta_i}{\partial u_i}
=\frac{\partial A_{ijk}}{\partial u_i}+\frac{\partial \theta_j}{\partial u_i}+\frac{\partial \theta_k}{\partial u_i}
=\frac{\partial \theta_j}{\partial u_i}\cosh l_{ij}+\frac{\partial \theta_k}{\partial u_i}\cosh l_{ik}.
\end{equation}
The formula (\ref{relation of partial derivative hyperbolic}) implies $\frac{\partial \theta_i}{\partial u_i}\rightarrow -\infty$
as $(f_i, f_j, f_k)\rightarrow (\overline{f}_i, \overline{f}_j, \overline{f}_k)\in \partial V_i$.
\end{remark}

Lemma \ref{hyperbolic symmetry} shows that the Jacobian matrix
\begin{equation*}
\begin{aligned}
\Lambda^H_{ijk}:=\frac{\partial (\theta_i, \theta_j, \theta_k)}{\partial ( u_i, u_j, u_k)}
=\left(
   \begin{array}{ccc}
     \frac{\partial \theta_i}{\partial u_i} & \frac{\partial \theta_i}{\partial u_j} & \frac{\partial \theta_i}{\partial u_k} \\
     \frac{\partial \theta_j}{\partial u_i} & \frac{\partial \theta_j}{\partial u_j} & \frac{\partial \theta_j}{\partial u_k} \\
     \frac{\partial \theta_k}{\partial u_i} & \frac{\partial \theta_k}{\partial u_j} & \frac{\partial \theta_k}{\partial u_k} \\
   \end{array}
 \right)
\end{aligned}
\end{equation*}
is symmetric.
Furthermore, we have the following result on the rank of the Jacobian matrix  $\Lambda^H_{ijk}$.

\begin{lemma}\label{hyperbolic rank lemma}
For the weighted triangle $(\{ijk\}, \varepsilon,\eta)$,
the rank of the Jacobian matrix $\Lambda^H_{ijk}$ is $3$
for any nondegenerate hyperbolic discrete conformal factor on $(\{ijk\}, \varepsilon,\eta)$.
\end{lemma}
\proof
By the chain rules, we have
\begin{equation}\label{hyperbolic chain rules proof for rank}
\begin{aligned}
\frac{\partial (\theta_i, \theta_j, \theta_k)}{\partial ( u_i, u_j, u_k)}
=\frac{\partial (\theta_i, \theta_j, \theta_k)}{\partial ( l_{jk}, l_{ik}, l_{ij})}
\cdot\frac{\partial ( l_{jk}, l_{ik}, l_{ij})}{\partial ( u_i, u_j, u_k)}.
\end{aligned}
\end{equation}
The derivative cosine law (\cite{CL}, Lemma A1) gives
\begin{equation*}
\begin{aligned}
\frac{\partial (\theta_i, \theta_j, \theta_k)}{\partial ( l_{jk}, l_{ik}, l_{ij})}
=\frac{1}{A}\left(
              \begin{array}{ccc}
                \sinh l_{jk} &   &   \\
                  & \sinh l_{ik} &   \\
                  &   & \sinh l_{ij} \\
              \end{array}
            \right)
            \left(
              \begin{array}{ccc}
                1 & -\cos\theta_k & -\cos\theta_j \\
                -\cos\theta_k & 1 & -\cos\theta_i \\
                -\cos\theta_j & -\cos\theta_i & 1 \\
              \end{array}
            \right).
\end{aligned}
\end{equation*}
This implies
\begin{equation*}
\begin{aligned}
&\det\left(\frac{\partial (\theta_i, \theta_j, \theta_k)}{\partial ( l_{jk}, l_{ik}, l_{ij})}\right)\\
=&\frac{\sinh l_{ij}\sinh l_{ik}\sinh l_{jk}}{A^3}
  \det \left(
              \begin{array}{ccc}
                1 & -\cos\theta_k & -\cos\theta_j \\
                -\cos\theta_k & 1 & -\cos\theta_i \\
                -\cos\theta_j & -\cos\theta_i & 1 \\
              \end{array}
            \right)\\
=&-\frac{\sinh l_{ij}\sinh l_{ik}\sinh l_{jk}}{A^3}(-1+\cos\theta_i^2+\cos\theta_j^2+\cos\theta_k^2+2\cos \theta_i\cos \theta_j\cos \theta_k)\\
=&-\frac{4\sinh l_{ij}\sinh l_{ik}\sinh l_{jk}}{A^3}\\
 &\cdot\cos\frac{\theta_i+\theta_j+\theta_k}{2}\cos\frac{\theta_i+\theta_j-\theta_k}{2}
   \cos\frac{\theta_i-\theta_j+\theta_k}{2}\cos\frac{\theta_i-\theta_j-\theta_k}{2}.
\end{aligned}
\end{equation*}
By the area formula for hyperbolic triangles, we have
$\theta_i+\theta_j+\theta_k\in (0,\pi)$.
This implies
$$\frac{\theta_i+\theta_j+\theta_k}{2}, \frac{\theta_i+\theta_j-\theta_k}{2},
\frac{\theta_i-\theta_j+\theta_k}{2}, \frac{\theta_i-\theta_j-\theta_k}{2}\in (-\frac{\pi}{2}, \frac{\pi}{2}).$$
Then we have
\begin{equation}\label{determinant of theta w.r.t l}
  \det (\frac{\partial (\theta_i, \theta_j, \theta_k)}{\partial ( l_{jk}, l_{ik}, l_{ij})})<0.
\end{equation}

By (\ref{defn of hyperbolic length}) and (\ref{derivative of fi in ui hyperbolic}), we have
\begin{equation*}
\begin{aligned}
&\frac{\partial ( l_{jk}, l_{ik}, l_{ij})}{\partial ( u_i, u_j, u_k)}\\
=&\left(
    \begin{array}{ccc}
      \frac{1}{\sinh l_{jk}} &   &   \\
        & \frac{1}{\sinh l_{ik}} &   \\
        &   & \frac{1}{\sinh l_{ij}} \\
    \end{array}
  \right)\\
 &\cdot
   \left(
     \begin{array}{ccc}
       0 & \varepsilon_jS_j^2C_k+\eta_{jk}S_jS_kC_j & \varepsilon_kS_k^2C_j+\eta_{jk}S_jS_kC_k \\
       \varepsilon_iS_i^2C_k+\eta_{ik}S_iS_kC_i & 0 & \varepsilon_kS_k^2C_i+\eta_{ik}S_iS_kC_k \\
       \varepsilon_iS_i^2C_j+\eta_{ij}S_iS_jC_i  & \varepsilon_jS_j^2C_i+\eta_{ij}S_iS_jC_j & 0 \\
     \end{array}
   \right).
\end{aligned}
\end{equation*}
This implies
\begin{equation}\label{Jac l of u positive}
\begin{aligned}
&\det\left(\frac{\partial ( l_{jk}, l_{ik}, l_{ij})}{\partial ( u_i, u_j, u_k)}\right)\\
=&\frac{S_iS_jS_k}{\sinh l_{ij}\sinh l_{ik}\sinh l_{jk}}\\
 &\cdot[2(\varepsilon_i\varepsilon_j\varepsilon_k+\eta_{ij}\eta_{ik}\eta_{jk})S_iS_jS_kC_iC_jC_k+\gamma_iS_iC_i(\varepsilon_kS_k^2C^2_j+\varepsilon_jS_j^2C_k^2)\\
  &+\gamma_jS_jC_j(\varepsilon_iS_i^2C_k^2+\varepsilon_kS_k^2C_i^2)
  +\gamma_kS_kC_k(\varepsilon_iS_i^2C_j^2+\varepsilon_jS_j^2C_i^2)]\\
\geq&\frac{2S_i^2S_j^2S_k^2C_iC_jC_k}{\sinh l_{ij}\sinh l_{ik}\sinh l_{jk}}
 [\varepsilon_i\varepsilon_j\varepsilon_k+\eta_{ij}\eta_{ik}\eta_{jk}+\gamma_i\varepsilon_j\varepsilon_k+\gamma_j\varepsilon_i\varepsilon_k
    +\gamma_k\varepsilon_i\varepsilon_j]\\
=&\frac{2S_i^2S_j^2S_k^2C_iC_jC_k}{\sinh l_{ij}\sinh l_{ik}\sinh l_{jk}}
 (\varepsilon_i\varepsilon_j+\eta_{ij})(\varepsilon_i\varepsilon_k+\eta_{ik})
   (\varepsilon_j\varepsilon_k+\eta_{jk})\\
>&0,
\end{aligned}
\end{equation}
where the structure condition (\ref{structure condition 2}) is used in the third line and
the structure condition (\ref{structure condition 1}) is used in the last line.

Therefore, by
(\ref{hyperbolic chain rules proof for rank}), (\ref{determinant of theta w.r.t l}) and (\ref{Jac l of u positive}),
we have $\det (\frac{\partial (\theta_i, \theta_j, \theta_k)}{\partial ( u_i, u_j, u_k)})<0$. This implies the rank of
the Jacobian matrix $\Lambda^H_{ijk}=\frac{\partial (\theta_i, \theta_j, \theta_k)}{\partial ( u_i, u_j, u_k)}$
is $3$.
\qed

As a consequence of Lemma \ref{hyperbolic symmetry} and Lemma \ref{hyperbolic rank lemma},
we have the following result on the negative definiteness of
the Jacobian matrix $\Lambda^H_{ijk}=\frac{\partial (\theta_i, \theta_j, \theta_k)}{\partial ( u_i, u_j, u_k)}$.
\begin{theorem}\label{hyperbolic Jacobian negativity}
For the weighted triangle $(\{ijk\}, \varepsilon,\eta)$,
the Jacobian matrix $\Lambda^H_{ijk}=\frac{\partial (\theta_i, \theta_j, \theta_k)}{\partial ( u_i, u_j, u_k)}$
is symmetric and negative definite for any nondegenerate
hyperbolic discrete conformal factor.
\end{theorem}
\proof
By Lemma \ref{hyperbolic rank lemma}, all the three eigenvalues of the Jacobian matrix $\Lambda^H_{ijk}$ are nonzero.
Taking $\Lambda^H_{ijk}$ as a matrix-valued function of $(f_i, f_j, f_k, \eta_{ij}, \eta_{ik}, \eta_{jk})\in \Omega^H_{ijk}$.
By the continuity of the eigenvalues of $\Lambda^H_{ijk}$ and the connectivity of the parameterized admissible space $\Omega^H_{ijk}$ in Corollary
\ref{coro connectivity of para admi space hyperbolic}, to prove the negative definiteness of
$\Lambda^H_{ijk}$,
we just need to find a point $p\in \Omega^H_{ijk}$ such that the eigenvalues of
$\Lambda^H_{ijk}$ at $p$ are negative.
Taking $p=(f_i, f_j, f_k, \eta_{ij}, \eta_{ik}, \eta_{jk})=(0, 0, 0, 1, 1, 1)$.
By Lemma \ref{hi sign at special point hyperbolic},
$p\in \Omega_{ijk}^H$ and $h_i(p)>0, h_j(p)>0, h_k(p)>0$.
Combining this with Lemma \ref{hyperbolic symmetry}, we have $\frac{\partial \theta_i}{\partial u_j}=\frac{\partial \theta_j}{\partial u_i}>0$
and
$\frac{\partial \theta_i}{\partial u_k}=\frac{\partial \theta_k}{\partial u_i}>0$ at $p$.
By (\ref{relation of partial derivative hyperbolic}), we have
$-\frac{\partial \theta_i}{\partial u_i}
>\frac{\partial \theta_j}{\partial u_i}+\frac{\partial \theta_k}{\partial u_i}$
at $p$.
By Lemma \ref{diagonal dominant} (a), the Jacobian matrix $\Lambda^H_{ijk}$ is negative definite and has three negative eigenvalues at $p$.
\qed

As a corollary of Theorem \ref{hyperbolic Jacobian negativity}, we have the following result on the
Jacobian matrix $\Lambda^H=\frac{\partial (K_1,\cdots,K_N)}{\partial (u_1,\cdots, u_N)}$
for nondegenerate hyperbolic discrete conformal structures.

\begin{corollary}\label{hyperbolic curvature Jacobian positivity}
Suppose $(M, \mathcal{T}, \varepsilon, \eta)$ is a weighted triangulated surface with
the weights $\varepsilon: V\rightarrow \{0, 1\}$ and $\eta: E\rightarrow \mathbb{R}$
satisfying the structure conditions (\ref{structure condition 1}) and (\ref{structure condition 2}).
Then the Jacobian matrix $\Lambda^H=\frac{\partial (K_1,\cdots,K_N)}{\partial (u_1,\cdots, u_N)}$
is symmetric and positive definite
for all nondegenerate hyperbolic discrete conformal factors on $(M, \mathcal{T}, \varepsilon, \eta)$.
\end{corollary}

The proof for Corollary \ref{hyperbolic curvature Jacobian positivity}
is the same as that for Corollary \ref{Euclidean curvature Jacobian positivity}, so
we omit the details of the proof here.

\subsection{Rigidity of hyperbolic discrete conformal structures}
Theorem \ref{hyperbolic admissible space} and Lemma \ref{hyperbolic symmetry} imply the following Ricci energy function for the weighted triangle $(\{ijk\}, \varepsilon,\eta)$
\begin{equation}\label{hyperbolic Ricci energy function for triangle}
\begin{aligned}
\mathcal{E}_{ijk}(u_i,u_j,u_k)=\int_{(\overline{u}_i, \overline{u}_j, \overline{u}_k)}^{(u_i, u_j, u_k)}\theta_idu_i+\theta_jdu_j+\theta_kdu_k
\end{aligned}
\end{equation}
is a well-defined smooth function on $\Omega^H_{ijk}(\eta)$ with $\nabla_{u_i}\mathcal{E}_{ijk}=\theta_i$.
The Ricci energy function $\mathcal{E}_{ijk}(u_i,u_j,u_k)$ was first constructed by Glickenstein-Thomas \cite{GT}
for Glickenstein's hyperbolic discrete conformal structures
under the assumption that the domain is simply connected.
Furthermore, Glickenstein-Thomas \cite{GT} used the Ricci energy function to prove a result on the local rigidity of Glickenstein's
hyperbolic discrete conformal structures.
For completeness, we give a sketch of  Glickenstein-Thomas's arguments here.
By Theorem \ref{hyperbolic Jacobian negativity}, $\mathcal{E}_{ijk}(u_i,u_j,u_k)$ is a locally strictly
concave function defined on $\Omega^H_{ijk}(\eta)$.
Set
\begin{equation}\label{hyperbolic Ricci energy function}
\begin{aligned}
\mathcal{E}(u_1, \cdots, u_N)=2\pi\sum_{i\in V}u_i-\sum_{\{ijk\}\in F}\mathcal{E}_{ijk}(u_i,u_j,u_k)
\end{aligned}
\end{equation}
to be the Ricci energy function defined on the admissible space $\Omega^H$ of nondegenerate hyperbolic discrete conformal factors
for $(M, \mathcal{T}, \varepsilon, \eta)$.
By Corollary \ref{hyperbolic curvature Jacobian positivity},
$\mathcal{E}(u_1, \cdots, u_N)$ is a locally strictly convex function on $\Omega^H$ with
$\nabla_{u_i} \mathcal{E}=K_i$. By Lemma \ref{injectivity of convex function},
the local rigidity of hyperbolic discrete conformal structures on $(M, \mathcal{T}, \varepsilon, \eta)$ follows.

To prove the global rigidity of hyperbolic discrete conformal structures,
we need to extend the inner angles in a hyperbolic triangle $\{ijk\}$
defined for nondegenerate hyperbolic discrete conformal factors
to be a globally defined function for $(f_i, f_j, f_k)\in \mathbb{R}^3$.
Parallelling to Lemma \ref{Euclidean extension}, we have the following extension for inner angles of hyperbolic triangles.

\begin{lemma}\label{hyperbolic extension}
For the weighted triangle $(\{ijk\}, \varepsilon,\eta)$,  the inner angles $\theta_i, \theta_j, \theta_k$
defined for nondegenerate hyperbolic discrete conformal factors
can be extended by constants to be continuous functions $\widetilde{\theta}_i, \widetilde{\theta}_j, \widetilde{\theta}_k$
defined for $(f_i, f_j, f_k)\in \mathbb{R}^3$ by setting
\begin{equation}\label{extension of theta_i hyperbolic}
\begin{aligned}
\widetilde{\theta}_i(f_i,f_j,f_k)=\left\{
                       \begin{array}{ll}
                         \theta_i, & \hbox{if $(f_i,f_j,f_k)\in \Omega^H_{ijk}(\eta)$;} \\
                         \pi, & \hbox{if $(f_i,f_j,f_k)\in V_i$;} \\
                         0, & \hbox{otherwise.}
                       \end{array}
                     \right.
\end{aligned}
\end{equation}
\end{lemma}
\proof
By Theorem \ref{hyperbolic admissible space},
$\Omega_{ijk}^H(\eta)=\mathbb{R}^3\setminus\sqcup_{\alpha\in \Lambda}V_\alpha$,
where
$\Lambda=\{q\in \{i,j,k\}|A_q=\eta_{st}^2-\varepsilon_s\varepsilon_t>0, \{q, s, t\}=\{i, j, k\}\}$ and
 $V_\alpha$ is a closed region in $\mathbb{R}^3$ bounded by the analytical function in (\ref{discription of V_i hyperbolic})
 defined on $\mathbb{R}^2$.

If $\Lambda= \varnothing$, then $\Omega_{ijk}^H(\eta)=\mathbb{R}^3$ and
$\theta_i, \theta_j, \theta_k$ is defined for all $(f_i, f_j, f_k)\in \mathbb{R}^3$.

If $\Lambda\neq \varnothing$, let $V_i$ be a connected component of $\mathbb{R}^3\setminus \Omega_{ijk}^H(\eta)$.
Suppose $(f_i, f_j, f_k)\in \Omega_{ijk}^H(\eta)$ tends to a point
$(\overline{f}_i, \overline{f}_j, \overline{f}_k)\in\partial V_i$ in $\mathbb{R}^3$.
By direct calculations, we have
\begin{equation}\label{expansion of hyperbolic triangle inequalities 2}
\begin{aligned}
&4\sinh\frac{l_{ij}+l_{ik}+l_{jk}}{2}\sinh\frac{l_{ij}+l_{ik}-l_{jk}}{2}\sinh\frac{l_{ij}-l_{ik}+l_{jk}}{2}\sinh\frac{-l_{ij}+l_{ik}+l_{jk}}{2}\\
=&(\cosh (l_{jk}+l_{ik})-\cosh l_{ij})(\cosh l_{ij}-\cosh (l_{jk}-l_{ik}))\\
=&(\cosh^2 l_{jk}-1)(\cosh l_{ik}^2-1)-(\cosh l_{jk}\cosh l_{ik}-\cosh l_{ij})^2\\
=&\sinh^2l_{jk}\sinh^2l_{ik} -\sinh^2l_{jk}\sinh^2l_{ik} \cos^2\theta_k\\
=&\sinh^2l_{jk}\sinh^2l_{ik}\sin^2\theta_k.
\end{aligned}
\end{equation}
Combining Lemma \ref{lemma hyper triangle ineq in Q}, (\ref{expansion of hyperbolic triangle inequalities})
and  the hyperbolic sine law,
(\ref{expansion of hyperbolic triangle inequalities 2}) implies that $\theta_i, \theta_j, \theta_k$ tends to $0$ or $\pi$ as
$(f_i, f_j, f_k)\rightarrow (\overline{f}_i, \overline{f}_j, \overline{f}_k)\in \partial V_i$.

By Remark \ref{hi positive in Vi hyperbolic}, for $(\overline{f}_i, \overline{f}_j, \overline{f}_k)\in \partial V_i$, we have
$h_i<0, h_j>0, h_k>0$.
By the continuity of $h_i, h_j, h_k$ of $(f_i,f_j,f_k)\in \mathbb{R}^3$,
there exists some neighborhood $U$
of $(\overline{f}_i, \overline{f}_j, \overline{f}_k)$ in $\mathbb{R}^3$
such that $h_i<0$, $h_j>0$, $h_k>0$
for $(f_i, f_j, f_k)\in \Omega^H_{ijk}(\eta)\cap U$.
This implies
$\frac{\partial \theta_j}{\partial f_i}=\frac{\partial \theta_j}{\partial u_i}\frac{1}{C_i}=\frac{S_i^2S_j^2S_k}{AC_i\sinh^2 l_{ij}}h_k>0$
for $(f_i,f_j,f_k)\in \Omega^H_{ijk}(\eta)\cap U$.
Similarly, $\frac{\partial \theta_k}{\partial f_i}>0$ for $(f_i,f_j,f_k)\in \Omega^H_{ijk}(\eta)\cap U$.
By the form (\ref{discription of V_i hyperbolic}) of $V_i$, we have $\theta_j, \theta_k\rightarrow 0$ as
$(f_i, f_j, f_k)\rightarrow (\overline{f}_i, \overline{f}_j, \overline{f}_k)\in \partial V_i$.
Otherwise, if  $\theta_j\rightarrow \pi$,  we have $\theta_j>\pi$ for
$(f_i+\epsilon, f_j, f_k)\in \Omega^H_{ijk}(\eta)\cap U$, $\epsilon>0$. It is impossible. The same arguments apply to $\theta_k$.

Furthermore, we have the following formula (\cite{V} page 66)
\begin{equation*}
\begin{aligned}
\tan^2\frac{A_{ijk}}{4}
=&\tanh \frac{l_{ij}+l_{ik}+l_{jk}}{2}\tanh \frac{l_{ij}+l_{ik}-l_{jk}}{2}\tanh \frac{l_{ij}-l_{ik}+l_{jk}}{2}\tanh \frac{-l_{ij}+l_{ik}+l_{jk}}{2}
\end{aligned}
\end{equation*}
for the area $A_{ijk}$ of the nondegenerate hyperbolic triangle $\{ijk\}$.  Combining this with the equation (\ref{expansion of hyperbolic triangle inequalities}), we have
\begin{equation}\label{hyperbolic area in QH}
\begin{aligned}
\tan^2\frac{A_{ijk}}{4}
=&\frac{\sinh\frac{l_{ij}+l_{ik}+l_{jk}}{2}\sinh\frac{l_{ij}+l_{ik}-l_{jk}}{2}\sinh\frac{l_{ij}-l_{ik}+l_{jk}}{2}\sinh\frac{-l_{ij}+l_{ik}+l_{jk}}{2}}
{16\cosh^2 \frac{l_{ij}+l_{ik}+l_{jk}}{4}\cosh^2 \frac{l_{ij}+l_{ik}-l_{jk}}{4}\cosh^2 \frac{l_{ij}-l_{ik}+l_{jk}}{4}\cosh^2 \frac{-l_{ij}+l_{ik}+l_{jk}}{4}}\\
=&\frac{S_i^2S_j^2S_k^2Q^H}{64\cosh^2 \frac{l_{ij}+l_{ik}+l_{jk}}{4}\cosh^2 \frac{l_{ij}+l_{ik}-l_{jk}}{4}\cosh^2 \frac{l_{ij}-l_{ik}+l_{jk}}{4}\cosh^2 \frac{-l_{ij}+l_{ik}+l_{jk}}{4}}.
\end{aligned}
\end{equation}
The equation (\ref{hyperbolic area in QH}) implies $A_{ijk}\rightarrow 0$
as $(f_i, f_j, f_k)\rightarrow (\overline{f}_i, \overline{f}_j, \overline{f}_k)\in \partial V_i$.
By $A_{ijk}=\pi-\theta_i-\theta_j-\theta_k$ and $\theta_j, \theta_k\rightarrow 0$, we have $\theta_i\rightarrow \pi$
as $(f_i, f_j, f_k)\rightarrow (\overline{f}_i, \overline{f}_j, \overline{f}_k)\in \partial V_i$.
Similar arguments apply for the other connected components of $\mathbb{R}^3\setminus \Omega_{ijk}^H(\eta)$.

Therefore, the extension (\ref{extension of theta_i hyperbolic}) defines a continuous extension of  the inner angle functions $\theta_i, \theta_j, \theta_k$
 on $\mathbb{R}^3$.
\qed

One can also use (\ref{relation of partial derivative hyperbolic}) to
prove $\theta_i\rightarrow \pi$ as $(f_i, f_j, f_k)\rightarrow (\overline{f}_i, \overline{f}_j, \overline{f}_k)\in \partial V_i$.

By Lemma \ref{hyperbolic extension}, we can extend the combinatorial curvature function $K$ defined on $\Omega^H$ to be defined for
all $f\in \mathbb{R}^N$ by setting
$\widetilde{K}_i=2\pi-\sum_{\{ijk\}\in F}\widetilde{\theta}_i$,
where $\widetilde{\theta}_i$ is the extension of $\theta_i$ defined  by (\ref{extension of theta_i hyperbolic}).

Taking $\widetilde{\theta}_i, \widetilde{\theta}_j, \widetilde{\theta}_k$ as functions of $(u_i, u_j, u_k)$.
Then the extensions $\widetilde{\theta}_i, \widetilde{\theta}_j, \widetilde{\theta}_k$ of $\theta_i, \theta_j, \theta_k$
 are continuous functions of
$(u_i, u_j, u_k)\in V_i\times V_j\times V_k$, where $V_q=\mathbb{R}$ if $\varepsilon_q=0$ and
$V_q=\mathbb{R}_{<0}=(-\infty, 0)$ if $\varepsilon_q=1$ for $q\in \{i,j,k\}$.
Combining this with Theorem \ref{Luo's convex extention}, the locally concave function
$\mathcal{E}_{ijk}$ defined by (\ref{hyperbolic Ricci energy function for triangle})
for nondegenerated $(u_i,u_j,u_k)$ on $(\{ijk\}, \varepsilon,\eta)$
can be extended to be a $C^1$ smooth concave function
\begin{equation}\label{extension of hyperbolic Ricci energy function for triangle}
\begin{aligned}
\widetilde{\mathcal{E}}_{ijk}(u_i,u_j,u_k)=\int_{(\overline{u}_i, \overline{u}_j, \overline{u}_k)}^{(u_i, u_j, u_k)}\widetilde{\theta}_idu_i+\widetilde{\theta}_jdu_j+\widetilde{\theta}_kdu_k
\end{aligned}
\end{equation}
defined for all $(u_i,u_j,u_k)\in V_i\times V_j\times V_k$ with $\nabla_{u_i}\widetilde{\mathcal{E}}_{ijk}=\widetilde{\theta}_i$.
As a result, the locally convex function $\mathcal{E}$ defined by (\ref{hyperbolic Ricci energy function}) for nondegenerate $u$ on $(M, \mathcal{T}, \varepsilon, \eta)$ can be extended to
be a $C^1$ smooth convex function
\begin{equation}\label{extended hyperbolic energy function}
\begin{aligned}
\widetilde{\mathcal{E}}(u_1, \cdots, u_N)=2\pi\sum_{i\in V}u_i-\sum_{\{ijk\}\in F}\widetilde{\mathcal{E}}_{ijk}(u_i,u_j,u_k)
\end{aligned}
\end{equation}
defined on $\mathbb{R}^{N_1}\times \mathbb{R}^{N_2}_{<0}$ with $\nabla_{u_i} \widetilde{\mathcal{E}}=\widetilde{K}_i=2\pi-\sum\widetilde{\theta}_i$,
where $N_1$ is the number of vertices $v_i$ in $V$ with $\varepsilon_i=0$ and $N_2=N-N_1$.

We have the following result on the rigidity of hyperbolic discrete conformal structures,
which is a generalization of Theorem \ref{main rigidity introduction} (b).

\begin{theorem}\label{Thm rigidity of HDCS context}
Suppose $(M, \mathcal{T}, \varepsilon, \eta)$ is a weighted triangulated surface with
the weights $\varepsilon: V\rightarrow \{0, 1\}$ and $\eta: E\rightarrow \mathbb{R}$
satisfying the structure conditions (\ref{structure condition 1}) and (\ref{structure condition 2}).
If there exists a nondegenrate hyperbolic discrete conformal factor $f_A\in \Omega^H$ and
a hyperbolic discrete conformal factor $f_B\in \mathbb{R}^N$ such that
$K(f_A)=\widetilde{K}(f_B)$. Then $f_A=f_B$.
\end{theorem}

\proof
Set
\begin{equation*}
\begin{aligned}
\mathcal{F}(t)=\widetilde{\mathcal{E}}((1-t)u_A+tu_B)=2\pi\sum_{i=1}^N[(1-t)u_{A,i}+tu_{B,i}]+\sum_{\{ijk\}\in F}\mathcal{F}_{ijk}(t),
\end{aligned}
\end{equation*}
where
$\mathcal{F}_{ijk}(t)=-\widetilde{\mathcal{E}}_{ijk}((1-t)u_A+tu_B)$.
Then $\mathcal{F}(t)$ is a $C^1$ smooth convex function for $t\in [0, 1]$ with $\mathcal{F}'(0)=\mathcal{F}'(1)$.
This implies that
$\mathcal{F}'(t)\equiv \mathcal{F}'(0)$ for all $t\in [0, 1]$.
Note that the admissible space $\Omega^H$ of nondegenerate hyperbolic discrete conformal factors is an
open subset of $\mathbb{R}^{N_1}\times \mathbb{R}^{N_2}_{<0}$, there exists $\epsilon>0$ such that $(1-t)u_A+tu_B$
corresponds to a nondegenerate hyperbolic discrete conformal factor for $t\in [0, \epsilon]$.
Note that $\mathcal{F}(t)$ is smooth for $t\in [0, \epsilon]$, we have
\begin{equation*}
\begin{aligned}
\mathcal{F}''(t)=(u_B-u_A) \Lambda^H (u_B-u_A)^T=0, \ \forall t\in [0, \epsilon].
\end{aligned}
\end{equation*}
By Corollary \ref{hyperbolic curvature Jacobian positivity}, this implies $u_A=u_B$.
As the transformation $u=u(f)$ defined by (\ref{ui in fi hyperbolic}) is a diffeomorphism, we have $f_A=f_B$.
\qed


\section{Relationships of Glickenstein's discrete conformal structures on surfaces and 3-dimensional hyperbolic geometry}\label{section 4}
\subsection{Construction of Glickenstein's discrete conformal structures via generalized hyperbolic tetrahedra}
The deep relationships of discrete conformal structures on polyhedral surfaces and $3$-dimensional hyperbolic geometry
were first discovered by Bobenko-Pinkall-Springborn \cite{BPS}
in the case of Luo's vertex scaling.  The relationships were then further studied in \cite{ZGZLYG}.
In this subsection, we study some more general cases.


We use the Klein model for $\mathbb{H}^3$ and take $\mathbb{S}^2$ as the ideal boundary $\partial \mathbb{H}^3$ of $\mathbb{H}^3$.
Suppose $\{ijk\}$ is a Euclidean or hyperbolic triangle generated by Glickenstein's discrete conformal structures in Definition \ref{defn of discrete conformal structure}.
The Ricci energy for the triangle $\{ijk\}$ is closely related to the co-volume of a generalized tetrahedron $T_{Oijk}$ in the
extended hyperbolic space $\mathbb{H}^3$, whose vertices are possibly truncated by a hyperbolic plane in $\mathbb{H}^3$.
In the following, we briefly describe the construction of $T_{Oijk}$ for $\varepsilon\in \{0,1\}$.
One can also refer to  \cite{BPS,ZGZLYG} for more information.

The generalized tetrahedron $T_{Oijk}$ has $4$ vertices $O, v_i, v_j, v_k$, which are ideal or hyper-ideal.
The vertex $O$ is called the bottom vertex.
\begin{description}
  \item[(1)] For the Euclidean background geometry, $O$ is ideal, i.e. $O\in \partial \mathbb{H}^3$. Please refer to Figure \ref{Euclidean-tetrahedron}.
  The Euclidean triangle $\{ijk\}$ is the intersection of the generalized hyperbolic tetrahedron $T_{Oijk}$
  with the horosphere $H_O$ at $O$.
  For the hyperbolic background geometry, $O$ is hyper-ideal, i.e. $O\not\in \mathbb{H}^3\cup \partial \mathbb{H}^3$,
  and the generalized hyperbolic tetrahedron is truncated by a hyperbolic plane $P_O$ in $\mathbb{H}^3$ dual to $O$.
  More precisely, $P_O=O^{\perp}\cap \mathbb{H}^3$, where $O^{\perp}$ is the timelike subspace in $\mathbb{R}^{3+1}$ orthogonal to the spacelike vector $O$.
  Please refer to Figure \ref{hyperbolic-tetrahedron}.
  The hyperbolic triangle $\{ijk\}$ is the intersection of the hyperbolic plane $P_O$ with the generalized hyperbolic tetrahedron $T_{Oijk}$.
\begin{figure}[!htb]
\begin{minipage}[t]{0.5\linewidth}
\centering
  \includegraphics[height=5cm,width=4cm]{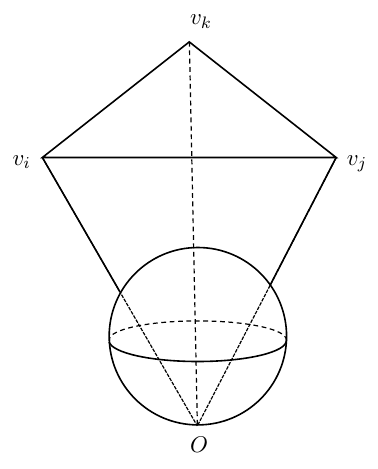}
  \caption{Tetrahedron for PL metric}
  \label{Euclidean-tetrahedron}
\end{minipage}
\hfill
\begin{minipage}[t]{0.5\linewidth}
\centering
  \includegraphics[height=5.2cm,width=4cm]{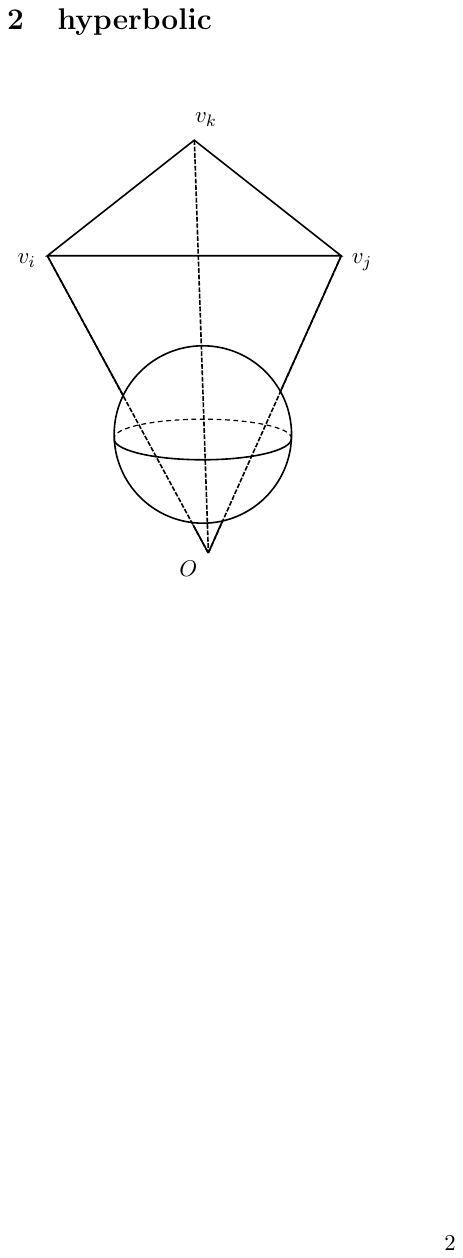}
  \caption{Tetrahedron for PH metric}
  \label{hyperbolic-tetrahedron}
\end{minipage}
\end{figure}

  \item[(2)] For $v_q\in \{v_i, v_j, v_k\}$, if the corresponding $\varepsilon_q=1$, then the vertex $v_q$ is hyper-ideal and
  the generalized tetrahedron $T_{Oijk}$ is truncated by a hyperbolic plane $P_q$ in $\mathbb{H}^3$ dual to $v_q$.
  If $O$ is also hyper-ideal, then it is required that $P_O\cap P_q=\emptyset$,
  which is equivalent to the line segment $Ov_q$ has nonempty intersection with $\mathbb{H}^3$ in the Klein model.
  If $\varepsilon_q=0$, then the vertex $v_q$ is ideal and we have a horosphere $H_q$ attached to $v_q$.
  For simplicity, we choose the horosphere $H_q$ so that it has no intersection
  with the hyperplane or horospheres attached to the other vertices of $T_{Oijk}$.
  \item[(3)] The signed edge length of $Ov_i, Ov_j, Ov_k$ are $-u_i, -u_j, -u_k$ respectively.
  \item[(4)] For the edge $v_iv_j$ in the extended hyperbolic space, the weight $\eta_{ij}$ is assigned as follows.
  \begin{description}
    \item[(a)] If $v_i, v_j$ are hyper-ideal and spans a spacelike or lightlike subspace $P_{ij}$, then $\eta_{ij}=\cos\beta_{ij}$,
    where $\beta_{ij}$ is determined by $-v_i\circ v_j=||v_i||\cdot||v_j||\cdot\cos\beta_{ij}$.
    Here we take $v_i,v_j$ as points in the Minkowski space, $\circ$ is the Lorentzian inner product in the Minkowski space
    and $||\cdot||$ is the norm of a spacelike vector.
    In fact, in the case that $v_i, v_j$ spans a spacelike subspace, the hyperbolic planes $P_i$ and $P_j$, dual to $v_i$ and $v_j$ respectively, intersect in $\mathbb{H}^3$
    and $\beta_{ij}$ is the dihedral angle
    determined by $P_i$ and $P_j$ in the truncated tetrahedron.
    \item[(b)] If $v_i, v_j$ are hyper-ideal and spans a timelike subspace, then $P_i$ and $P_j$ do not intersect in $\mathbb{H}^3$.
    Denote $\lambda_{ij}$ as the hyperbolic distance of $P_i$ and $P_j$, then $\eta_{ij}=\cosh \lambda_{ij}=-\frac{v_i\circ v_j}{||v_i||\cdot||v_j||}$.
    \item[(c)] If $v_i, v_j$ are ideal, we choose the horospheres $H_i, H_j$ at $v_i, v_j$ with $H_i\cap H_j=\emptyset$
    and set $\lambda_{ij}$ to be the distance from $H_i\cap \overline{v_iv_j}$ to $H_j\cap \overline{v_iv_j}$, where
    $\overline{v_iv_j}$ is the geodesic from $v_i$ to $v_j$. Then $\eta_{ij}=\frac{1}{2}e^{\lambda_{ij}}$.
    \item[(d)] If $v_i$ is ideal and $v_j$ is hyper-ideal, we choose the horosphere $H_i$ at $v_i$ to have no intersection with the
    hyperbolic plan $P_j$ dual to $v_j$. Set $\lambda_{ij}$ to be the distance from $H_i$ to $P_j$. Then $\eta_{ij}=\frac{1}{2}e^{\lambda_{ij}}$.
  \end{description}
\end{description}

In this setting, it can be checked that
the lengths for the edges in the Euclidean triangle $H_O\cap T_{Oijk}$ and in the hyperbolic triangle
$P_O\cap T_{Oijk}$ are given by (\ref{defn of Euclidean length}) and (\ref{defn of hyperbolic length}) respectively, where
$u_i=f_i$ for the Euclidean background geometry and $u_i$ is defined by (\ref{ui in fi hyperbolic}) in terms of $f_i$
for the hyperbolic background geometry.

By the hyperbolic cosine laws for generalized hyperbolic triangle $\{v_iv_jv_k\}$, it can be checked
that
$\eta_{st}+\varepsilon_s\varepsilon_t>0,\ \varepsilon_s\eta_{tq}+\eta_{st}\eta_{sq}>0, \ \{s,t,q\}=\{i,j,k\}.$
This proves Theorem \ref{theorem convexity of co-volume} (a).
We suggest the readers to refer to Appendix A of \cite{GL} for a full list of
formulas of hyperbolic sine and cosine laws for generalized hyperbolic triangles
used here.
In the case that $\varepsilon_i=\varepsilon_j=\varepsilon_k=1$, this can be proved in a geometric approach.
Note that $\eta_{ij}=-\frac{v_i\circ v_j}{||v_i||\cdot||v_j||}$ in this case.
Taking $\eta_{ij}+\eta_{ik}\eta_{jk}>0$ for example. Note that
\begin{equation}\label{Lorentz cross product}
\begin{aligned}
\frac{(v_j\otimes v_k)\circ (v_k\otimes v_i)}{||v_i||\cdot||v_j||\cdot||v_k||^2}
=\left(\frac{v_j}{||v_j||}\otimes \frac{v_k}{||v_k||}\right)\circ \left(\frac{v_k}{||v_k||}\otimes \frac{v_i}{||v_i}\right)
=-(\eta_{ij}+\eta_{ik}\eta_{jk}),
\end{aligned}
\end{equation}
where $\otimes$ is the Lorentzian cross product defined by $x\otimes y=J(x\times y)$ with $J=\text{diag}\{-1, 1,1\}$ for $x,y\in \mathbb{R}^3$.
Please refer to \cite{Rat} (Chapter 3) for more details on Lorentzian cross product.
By (\ref{Lorentz cross product}), to prove $\eta_{ij}+\eta_{ik}\eta_{jk}>0$, we just need to prove $(v_j\otimes v_k)\circ (v_k\otimes v_i)<0$.
In the following,
we use $P_{st}=\text{span}(v_s,v_t)$ to denote the two dimensional plane spanned by $v_s$ and $v_t$ in the Minkowski space.
By symmetry, we just need to consider the following six cases.

  \begin{figure}[!htb]
\centering
  \includegraphics[height=.5\textwidth,width=0.7\textwidth]{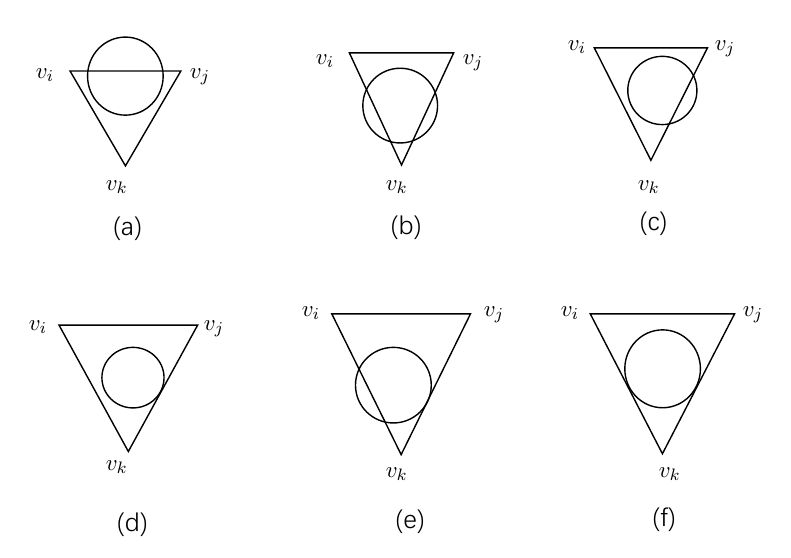}
  \caption{Generalized triangles}
  \label{Generalized triangles}
\end{figure}

\begin{description}
  \item[(a)] If $P_{ik}$ and $P_{jk}$ are spacelike, then $v_j\otimes v_k$, $v_k\otimes v_i$ are timelike with the same parity. This implies $(v_j\otimes v_k)\circ (v_k\otimes v_i)<0$. Please refer to Figure \ref{Generalized triangles} (a).
  \item[(b)] If $P_{ik}$ and $P_{jk}$ are timelike, then $v_j\otimes v_k$, $v_k\otimes v_i$ are spacelike and
  $(v_j\otimes v_k)\circ (v_k\otimes v_i)=-||v_j\otimes v_k||\cdot ||v_k\otimes v_i||\cdot\cosh d(P_{ik}, P_{jk})<0$.
  Please refer to Figure \ref{Generalized triangles} (b).
  \item[(c)] If $P_{ik}$  is spacelike and $P_{jk}$ is timelike, then $v_j\otimes v_k$ is spacelike and $v_k\otimes v_i$ is timelike.
  Then $(v_j\otimes v_k)\circ (v_k\otimes v_i)=-||v_j\otimes v_k||\cdot |||v_k\otimes v_i|||\cdot\sinh d<0$, where
   $|||x|||$ denotes the absolute value of $||x||=(x\circ x)^{1/2}$ for a timelike vector $x$ and
   $d$ is the distance of $\frac{v_k\otimes v_i}{|||v_k\otimes v_i|||}$ to $P_{jk}$. Please refer to Figure \ref{Generalized triangles} (c).
  \item[(d)] If $P_{ik}$  is spacelike and $P_{jk}$ is lightlike, then $v_j\otimes v_k$ is lightlike and
  $v_k\otimes v_i$ is timelike with the same parity. This implies $(v_j\otimes v_k)\circ (v_k\otimes v_i)<0$. Please refer to Figure \ref{Generalized triangles} (d).
  \item[(e)] If $P_{ik}$ is timelike and $P_{jk}$ is lightlike,
             then $v_j\otimes v_k$ is lightlike, $v_k\otimes v_i$ is spacelike with the same parity as $v_j\otimes v_k$.
             Then $(v_j\otimes v_k)\circ (v_k\otimes v_i)<0$. Please refer to Figure \ref{Generalized triangles} (e).
  \item[(f)] If $P_{ik}$ is lightlike and $P_{jk}$ is lightlike,
             then $v_j\otimes v_k$ is lightlike and $v_k\otimes v_i$ is lightlike with the same parity as $v_j\otimes v_k$.
             Furthermore, $v_j\otimes v_k$ and $v_k\otimes v_i$ are linearly independent.
             Then $(v_j\otimes v_k)\circ (v_k\otimes v_i)<0$. Please refer to Figure \ref{Generalized triangles} (f).
\end{description}


\subsection{Convexities of co-volume of generalized hyperbolic tetrahedra}
For the generalized hyperbolic tetrahedron $T_{Oijk}$, we have attached it with a generalized hyperbolic polyhedron $\widetilde{P}$ by truncating it by the hyperbolic planes dual to the vertices $O, v_i, v_j, v_k$.
If $\widetilde{P}$ is a generalized hyperbolic polyhedron in $\mathbb{H}^3$ with finite volume, we set $\widetilde{P}=P$.
Otherwise,  the generalized hyperbolic polyhedron $\widetilde{P}$ has hyper-ideal vertices
and we need to further truncate
$\widetilde{P}$ to get a generalized hyperbolic polyhedron $P$ with finite volume.
For example, in the case that $\varepsilon_i=\varepsilon_j=\varepsilon_k=1$, $P_{ij}, P_{ik}, P_{jk}$ are spacelike,
the generalized hyperbolic triangle $\triangle v_iv_jv_k$ has no intersection with $\partial \mathbb{H}^3$ and
the point $P_{ij}\cap P_{ik}\cap P_{jk}$ is hyper-ideal, we need to use a hyperbolic plane $P_{ijk}$ dual
to $P_{ij}\cap P_{ik}\cap P_{jk}$ to truncate $\widetilde{P}$
to get a finite hyperbolic polyhedron $P$.

Denote the volume of the generalized hyperbolic polyhedron $P$ by $V$.
By the Schl\"{a}fli formula \cite{R1}, we have
\begin{equation*}
\begin{aligned}
dV=-\frac{1}{2}(-u_id\theta_i-u_jd\theta_j-u_kd\theta_k+\lambda_{ij}d\beta_{ij}+\lambda_{ik}d\beta_{ik}+\lambda_{jk}d\beta_{jk}).
\end{aligned}
\end{equation*}
If $v_q, v_s\in \{v_i, v_j, v_k\}$ are spacelike and $P_{qs}$ is non-timelike, then $\beta_{qs}$ is fixed, otherwise $\lambda_{qs}$ is fixed.
Set
\begin{equation*}
\begin{aligned}
\mu_{qs}=\left\{
           \begin{array}{ll}
             0, & \hbox{ if $\varepsilon_q=\varepsilon_s=1$ and $P_{qs}$ is non-timelike,} \\
             1, & \hbox{otherwise.}
           \end{array}
         \right.
\end{aligned}
\end{equation*}
Define the co-volume by
\begin{equation}\label{definition of covolume}
\begin{aligned}
\widehat{V}=2V-u_i\theta_i-u_j\theta_j-u_k\theta_k+\mu_{ij}\lambda_{ij}\beta_{ij}
+\mu_{ik}\lambda_{ik}\beta_{ik}+\mu_{jk}\lambda_{jk}\beta_{jk}.
\end{aligned}
\end{equation}
Then we have
\begin{equation}\label{differential of covolume}
\begin{aligned}
d\widehat{V}=-\theta_idu_i-\theta_jdu_j-\theta_kdu_k.
\end{aligned}
\end{equation}
By Theorem \ref{Euclidean Jacobian negativity} and Theorem \ref{hyperbolic Jacobian negativity},
the equation (\ref{differential of covolume}) implies the co-volume function $\widehat{V}$ is convex in $u_i, u_j, u_k$.
As a result, the co-volume function $\widehat{V}$ is convex in the edge lengths $l_{Ov_i}=-u_i, l_{Ov_j}=-u_j, l_{Ov_k}=-u_k$.
This completes the proof of Theorem \ref{theorem convexity of co-volume} (b).

\section{Open problems}\label{section 5}
\subsection{Convergence of Glickenstein's discrete conformal structures}
In \cite{Th2}, Thurston conjectured that the tangential circle packing can be used to approximate the Riemann mapping.
Thurston's conjecture was then proved elegantly by Rodin-Sullivan \cite{RS}.
Since then, there have been lots of important works on Thurston's conjecture. See \cite{H,HZ1,HZ2} and others.
For Luo's vertex scalings, the corresponding convergence to the Riemann mapping was recently proved by Luo-Sun-Wu \cite{LSW}.
See also \cite{GLW, L6, WZ} for related works.
For Bowers-Stephenson's inversive distance circle packings,
the corresponding convergence to Riemann mapping was recently proved by Chen-Luo-Xu-Zhang \cite{CLXZ}.
Note that Glickenstein's discrete conformal structures are natural generalizations
of Bowers-Stephenson's inversive circle packings and Luo's vertex scalings.
It is convinced that Thurston's conjecture is still true for Glickenstein's discrete conformal structures.

\subsection{Discrete uniformization theorems for Glickenstein's discrete conformal structures}
An interesting question about Glickenstein's discrete conformal structures on polyhedral surfaces is the existence of a discrete conformal factor with the prescribed combinatorial curvature.
In the special case that the prescribed combinatorial curvature is constant, this corresponds to
the discrete uniformization theorem.
For Luo's vertex scaling, the discrete uniformization theorems were established in \cite{GGLSW,GLSW,Sp2}.
Note that Luo's vertex scalings  correspond to $\varepsilon\equiv1$ in Glickenstein's discrete conformal structures.
This motivates us to study the discrete uniformization theorem for Glickenstein's discrete conformal structures.

Suppose $(M, V)$ is a marked surface and $V$ is a nonempty finite subset of $M$.
A weight $\varepsilon: V\rightarrow \{0,1\}$ is defined on $V$. The triple $(M, V, \varepsilon)$ is called a weighted marked surface.
Motivated by Glickenstein's works \cite{G1,G2a,G3,G4}, we have the following definition of weighted Delaunay triangulation.

\begin{definition}
Suppose $(M, V, \varepsilon)$ is a weighted marked surface with a PL metric $d$, and
$\mathcal{T}$ is a geometric triangulation of $(M, V, \varepsilon)$ with every triangle $\{ijk\}\in \mathcal{T}$ have a
well-defined geometric center $C_{ijk}$.
Suppose $\{ij\}$ is an edge shared by two adjacent Euclidean triangles $\{ijk\}$ and $\{ijl\}$.
The edge $\{ij\}$ is called weighted Delaunay if $h_{ij,k}+h_{ij,l}\geq 0$, where
$h_{ij,k}, h_{ij,l}$ are the signed distance of $C_{ijk}, C_{ijl}$ to the edge $\{ij\}$ respectively.
The triangulation $\mathcal{T}$ is called weighted Delaunay in $d$ if every edge in the triangulation is weighted Delaunay.
\end{definition}

One can also define the weighted Delaunay triangulation using the power distance in Remark \ref{power distance}.
For a PL metric $d$ on $(M, V, \varepsilon)$, its weighted Voronoi decomposition is defined to be the
connection of 2-cells $\{R(v)|v\in V\}$, where $R(v)=\{x\in M|\pi_v(x)\leq \pi_{v'}(x) \ \text{for all } v'\in V\}$ is defined by the
power distance.
The dual cell-decomposition $\mathcal{C}(d)$ of the weighted Voronoi decomposition
is called the weighted Delaunay tessellation of $(M, V, \varepsilon, d)$.
A weighted Delaunay triangulation $\mathcal{T}$ of $(M, V, \varepsilon, d)$ is  a geometric
triangulation of the weighted Delaunay tessellation $\mathcal{C}(d)$ by further triangulating
all non-triangular 2-dimensional cells without introducing extra vertices.
As the power distance is a generalization of Euclidean distance, the weighted Delaunay triangulation is a generalization
of the Delaunay triangulation.

Following Gu-Luo-Sun-Wu \cite{GLSW}, we introduce the following new definition of discrete conformality, which allows the
triangulation of the weighted marked surface $(M, V, \varepsilon)$ to be changed.

\begin{definition}\label{new definition of edcf}
Two piecewise linear metrics $d, d'$ on $(M, V, \varepsilon)$ are discrete conformal if
there exist sequences of PL metrics $d_1=d$, $\cdots$,  $d_m=d'$
on $(M, V, \varepsilon)$ and triangulations $\mathcal{T}_1, \cdots, \mathcal{T}_m$ of
$(M, V, \varepsilon)$ satisfying
\begin{description}
  \item[(a)] (Weighted Delaunay condition) each $\mathcal{T}_i$ is weighted Delaunay in $d_i$,
  \item[(b)] (Discrete conformal condition) if $\mathcal{T}_i=\mathcal{T}_{i+1}$, there exists two functions
  $u_i, u_{i+1}:V\rightarrow \mathbb{R}$ such that if $e$ is an edge in $\mathcal{T}_i$ with end points $v$ and $v'$,
  then the lengths $l_{d_{i+1}}(e)$ and  $l_{d_{i}}(e)$ of $e$ in $d_i$ and $d_{i+1}$ are defined by (\ref{defn of Euclidean length})
  using $u_i$ and $u_{i+1}$ respectively with the same weight $\eta: E\rightarrow \mathbb{R}$.
  \item[(c)] if $\mathcal{T}_i\neq \mathcal{T}_{i+1}$, then $(S, d_i)$ is isometric to $(S, d_{i+1})$
  by an isometry homotopic to identity in $(S, V)$.
\end{description}
\end{definition}
The space of PL metrics on $(M, V, \varepsilon)$ discrete conformal to $d$ is called the conformal class of $d$ and
denoted by $\mathcal{D}(d)$.

Motivated by Gu-Luo-Sun-Wu's discrete uniformization theorem for vertex scalings of PL metrics in \cite{GLSW},
we have the following conjecture on
the discrete uniformization for Glickenstein's Euclidean discrete conformal structures on weighted marked surfaces.

\begin{conjecture}
Suppose $(M, V, \varepsilon)$ is a closed connected weighted marked surface with $\varepsilon: V\rightarrow \{0,1\}$, $\chi(M)=0$
and $d$ is a PL metric on $(M, V, \varepsilon)$.
There exists a PL metric $d'\in \mathcal{D}(d)$, unique up to scaling and isometry homotopic to the identity
 on $(M, V, \varepsilon)$, such that $d'$ is discrete conformal to $d$ and the combinatorial curvature of $d'$ is $0$.
\end{conjecture}

For the hyperbolic background geometry, one can define the corresponding weighted Delaunay triangulation and
the corresponding discrete conformality similarly.
We have the following conjecture on the discrete uniformization  for Glickenstein's hyperbolic discrete conformal structures on weighted marked surfaces.

\begin{conjecture}
Suppose $(M, V, \varepsilon)$ is a closed connected weighted marked surface with $\varepsilon: V\rightarrow \{0,1\}$,  $\chi(M)<0$ and $d$ is a PH metric on $(M, V, \varepsilon)$.
There exists a unique
PH metric $d'\in \mathcal{D}(d)$ on $(M, V,\varepsilon)$ so that $d'$ is discrete conformal
to $d$ and the combinatorial curvature of $d'$ is 0.
\end{conjecture}

One can also study the prescribing combinatorial curvature problem for Glickenstein's discrete conformal structures on polyhedral surfaces.
Results similar to the results in \cite{GGLSW, GLSW, T1} are convinced to be true
for Glickenstein's discrete conformal structures on polyhedral surfaces.

\end{document}